\documentclass[11pt, a4paper]{amsart}
\usepackage{preamble}

\begin{document}
\bibliographystyle{plain}
\title{Energy minimizing harmonic $2$-spheres in metric spaces}

\author{Damaris Meier}
\address
  {Department of Mathematics\\ University of Fribourg\\ Chemin du Mus\'ee 23\\ 1700 Fribourg, Switzerland}
\email{damaris.meier@unifr.ch}

\author{Noa Vikman}
\address
  {Department of Mathematics\\ University of Fribourg\\ Chemin du Mus\'ee 23\\ 1700 Fribourg, Switzerland}
\email{noa.vikman@unifr.ch}

\author{Stefan Wenger}
\address
  {Department of Mathematics\\ University of Fribourg\\ Chemin du Mus\'ee 23\\ 1700 Fribourg, Switzerland}
\email{stefan.wenger@unifr.ch}

\date{\today}

\thanks{Research supported by Swiss National Science Foundation Grant 212867 and partially funded by the Deutsche Forschungsgemeinschaft (DFG, German Research Foundation) under Germany's Excellence Strategy – EXC-2047/1 – 390685813.}

\begin{abstract}
   In their seminal 1981 article, Sacks-Uhlenbeck famously proved the existence of non-trivial harmonic $2$-spheres in every closed Riemannian manifold with non-zero second homotopy group. Their arguments heavily rely on PDE techniques. The purpose of the present paper is to develop a conceptually simple metric approach to the existence of harmonic spheres. This allows us to generalize the Sacks-Uhlenbeck result to a large class of compact metric spaces.
\end{abstract}

\maketitle

\section{Introduction} \label{sec:Introduction}
\subsection{Background}\label{sec:background}
The classical existence problem for harmonic maps asks whether a given continuous map between Riemannian manifolds $M$ and $N$ can be deformed by a homotopy to a harmonic map, and thus to a critical point of the Dirichlet energy functional. This and related problems have attracted a lot of research activity at least since the 1960s. We refer to \cite{EL78} and \cite{EL88} for early overviews of the topic. The existence problem has obtained particular attention in the case when $M$ is a closed $2$-dimensional Riemannian manifold, in the sequel referred to as a closed surface. Lemaire \cite{Lem78}, Schoen-Yau \cite{SY79}, Sacks-Uhlenbeck \cite{SU-1981} proved independently that if $N$ is a closed Riemannian manifold with trivial second homotopy group, then every homotopy class of maps from a closed surface $M$ to $N$ contains an energy minimizing harmonic map. As was shown in \cite{Futaki-1980}, \cite{SU-1981}, this fails in general without the condition on the second homotopy group. Nevertheless, Sacks-Uhlenbeck proved in their influential paper \cite{SU-1981} that if $N$ has non-trivial second homotopy group then there exists a generating set $P$ of $\pi_2(N)$ such that each element of $P$ contains an energy minimizing harmonic map in its free homotopy class; moreover, every such minimizer is a conformal branched immersion. This result, whose proof heavily relies on PDE techniques, has led to a flurry of research around the harmonic map problem in the past four decades.

In recent years, there has been an increasing interest in studying harmonic maps with values in a metric space. Korevaar-Schoen introduced a notion of Sobolev maps from Riemannian manifolds to complete metric spaces and proved the existence of harmonic maps into metric spaces of non-positive curvature (locally ${\rm CAT}(0)$ spaces) in \cite{KS93}; see also Jost \cite{Jost94} for related results around the same time. By now, there exists a robust theory of Sobolev maps from Riemannian and more general domains to complete metric spaces \cite{KS93}, \cite{Jost94}, \cite{Haj96}, \cite{HKST15}.

In the present paper we study the existence problem for harmonic $2$-spheres in the very general context of metric spaces admitting a local quadratic isoperimetric inequality for curves. 

\bd
    A complete metric space $X$ is said to admit a local quadratic isoperimetric inequality if there exist $C, l_0>0$ such that every Lipschitz curve $\gamma\colon S^1\to X$ of length $\length(\gamma)\leq l_0$ is the trace of a Sobolev map $u\in W^{1,2}(D, X)$ with $$\Area(u)\leq C\cdot \length(\gamma)^2,$$ where $D$ denotes the Euclidean unit disc.
\ed

We refer to \Cref{sec:Sobolev} for notions related to Sobolev maps to metric spaces. Many interesting classes of metric spaces admit a local quadratic isoperimetric inequality, see the examples listed after the statement of our main result, \Cref{thm:main-introduction}. Metric spaces with a local quadratic isoperimetric inequality provide a suitable context for studying area and energy minimization problems, see for example the recent papers \cite{LW15-Plateau}, \cite{LW17-en-area}, \cite{LW16-harmonic}, \cite{LW-intrinsic}. This has had applications in various areas, see for example \cite{LW-isoperimetric}, \cite{LW-param}, \cite{LS-2019}, \cite{LWY20}, \cite{LS-2020}.

The purpose of the present paper is to develop a metric approach to the existence problem and prove a strengthening of the Sacks-Uhlenbeck theorem and of related results in a large class of compact metric spaces admitting a local quadratic isoperimetric inequality. Our approach also yields a conceptually simple proof in the setting of smooth Riemannian manifolds. Furthermore, our main theorem and its consequences generalize and strengthen results in the smooth setting \cite{SU-1981}, \cite{Lem78}, \cite{SY79}, \cite{Jost-1991-two-dim-var-problems} as well as more recent results in the setting of compact locally ${\rm CAT}(1)$ spaces \cite{Breiner-et-al}.

\subsection{Statements of main results}\label{sec:intro-main-results}
Let $X$ be a complete metric space and $\varphi\colon M\to X$ a continuous map from  a closed surface $M$ (i.e.~a smooth closed $2$-dimensional Riemannian manifold) to $X$. Denote by $e(\varphi)$ the infimal energy of continuous Sobolev maps homotopic to $\varphi$. There are several natural notions of energy for metric space valued Sobolev maps, see \Cref{sec:energy-new}, and our results apply to all of them.  Throughout this introduction, we fix a definition of energy. Notice that the Korevaar-Schoen energy from \cite{KS93} corresponds to the Dirichlet energy when the ambient space $X$ is a Riemannian manifold. Under the assumptions of our main theorem, the homotopy class of $\varphi$ contains a representative of finite energy, see \Cref{prop:Sobolev-homotopic-to-cont}, but it need not contain an energy minimizer, see \Cref{example:dumbbell}. However, our main theorem will show that $\varphi$ has an \textit{iterated decomposition} into a finite number of continuous maps each of which admits an energy minimizer in its homotopy class. For the following definition, recall that a continuous map is said to be essential if it is not null-homotopic.

\bd
A pair of continuous maps $\varphi_0\colon M\to X$ and $\varphi_1\colon S^2\to X$ is said to decompose a continuous map $\varphi\colon M\to X$ if $\varphi_0$ agrees with $\varphi$ on the complement of a disc $B\subset M$, and $\varphi_1$ is essential and coincides with the map obtained from gluing the restrictions $\varphi|_{\bar{B}}$ and $\varphi_0|_{\bar{B}}$ along $\partial B$. If $M$ is diffeomorphic to $S^2$ then we also require that $\varphi_0$ is essential.
\ed

\begin{figure}
    \centering
    \includegraphics[width=\linewidth]{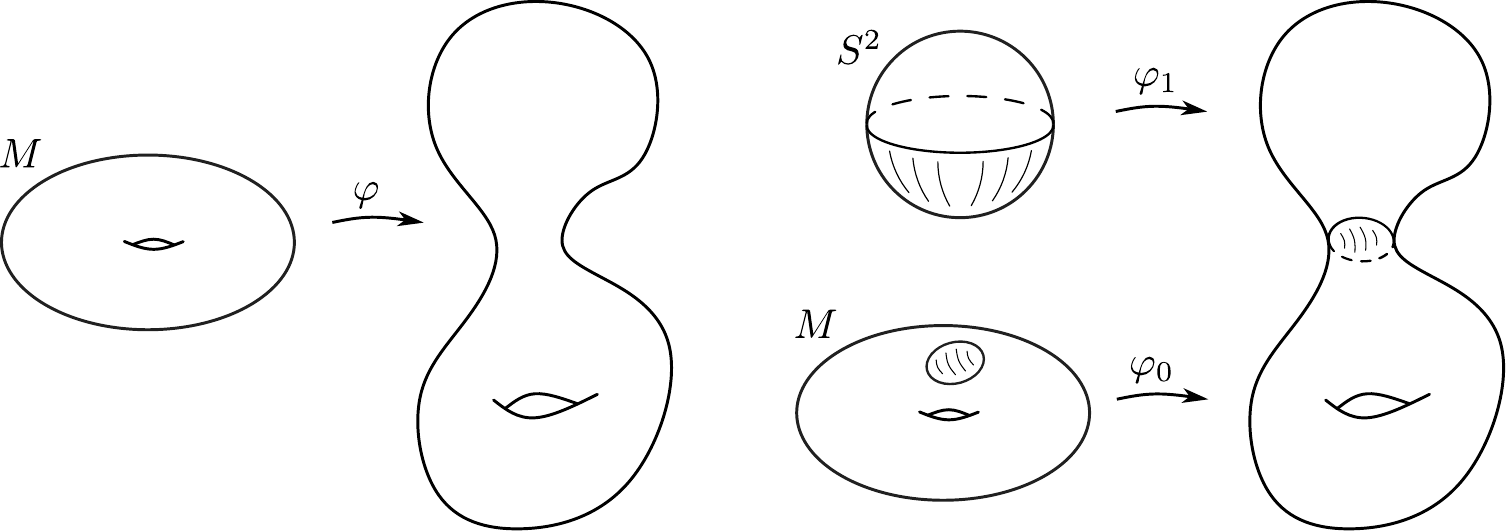}
    \caption{Decomposition $\varphi_0$ and $\varphi_1$ of a map $\varphi\colon M\to X$.}
    \label{fig:decomposition}
\end{figure}

We refer to \Cref{fig:decomposition} for an illustration of this concept. The process of repeatedly decomposing a map results in what we call an \textit{iterated decomposition}, inductively defined as follows. We say that $\varphi_0\colon M\to X$ is a zero-step iterated decomposition of $\varphi$ if $\varphi_0=\varphi$. Furthermore, for $k>0$, we define a $k$-step iterated decomposition of $\varphi$ as a $(k+1)$-tuple of maps $\varphi_0\colon M\to X$ and $\varphi_i\colon S^2\to X$, where $i=1,\dots,k$, obtained from a $(k-1)$-step iterated decomposition of $\varphi$ by taking a decomposition of one of its maps and keeping the remaining $k-1$ maps the same. It can be shown that in the setting of our main result the total energy $e(\varphi_0)  + \dots + e(\varphi_k)$ of an iterated decomposition is always at least $e(\varphi)$, see \Cref{lem:energy-of-decomposition}. 

We turn to our main result and first state the assumptions on our metric space. These will be in place throughout the introduction. Let $X$ be a compact and quasiconvex metric space admitting a local quadratic isoperimetric inequality. Assume furthermore that every continuous map from $S^2$ to $X$ of sufficiently small diameter is null-homotopic. 

\bt \label{thm:main-introduction}
Let $X$ be as above, and let $M$ be a closed surface, equipped with a Riemannian metric. Then every continuous map $\varphi\colon M\to X$ has an iterated decomposition satisfying 
\begin{equation}\label{eq:energy-quantization}
 e(\varphi_0) + e(\varphi_1) + \dots + e(\varphi_k) = e(\varphi)   
\end{equation}
and such that every $\varphi_i$ contains an energy minimizer in its homotopy class.
\et

Examples of spaces satisfying the assumptions of our theorem include closed Riemannian manifolds, more generally, compact Lipschitz manifolds, compact locally ${\rm CAT}(\kappa)$ spaces, compact Alexandrov spaces, some compact sub-Riemannian manifolds such as those locally modeled on higher Heisenberg groups, and many more. In particular, our theorem  generalizes and strengthens results in Sacks-Uhlenbeck \cite{SU-1981}, Jost \cite{Jost-1991-two-dim-var-problems}, and Breiner et al.~\cite{Breiner-et-al}. The equality \eqref{eq:energy-quantization} is also known as energy identity or energy quantization.

We now state several consequences of our main theorem. Since iterated decompositions are related to sums in the second homotopy group $\pi_2(X)$ we first obtain the following generalization to metric spaces of \cite[Theorem 5.9]{SU-1981}.

\bc\label{thm:full-Sacks-Uhlenbeck-metric-space}
 Let $X$ be as above. Then there exists a family $\Gamma$ of based continuous maps from $S^2$ to $X$ that generates $\pi_2(X)$ and such that every element of $\Gamma$ contains an energy minimizer in its (free) homotopy class.
\ec

When $\pi_2(X)$ is trivial then the only iterated decomposition of a map is the map itself. Therefore, as a direct consequence of \Cref{thm:main-introduction}, we get the following result which generalizes \cite[Theorem 5.8]{SU-1981}.

\bc\label{cor:trivial-second-homotopy-group-intro}
 Let $X$ be as above, and suppose furthermore that $\pi_2(X)$ is trivial. Let $M$ be a closed surface, equipped with a Riemannian metric. Then every continuous map from $M$ to $X$ contains an energy minimizer in its homotopy class.
\ec

We mention that this result can also be proved using the methods developed in \cite{Soultanis-Wenger-2022}, see \Cref{sec:1-homotopy-classes} below for a second proof along these lines.

We next discuss regularity properties of energy minimizers in homotopy classes. As is well-known, harmonic maps from a $2$-dimensional sphere into a Riemannian manifold are weakly conformal, see for example \cite{SU-1981}. Our theorem below shows that in the setting of metric spaces, energy minimizers in homotopy classes satisfy a metric variant of conformality called \textit{infinitesimal quasiconformality}. This means, roughly speaking, that infinitesimal balls get mapped to sets of bounded excentricity. We refer to \Cref{sec:energy-new} for the precise definition of infinitesimal quasiconformality.

\bt \label{thm:regularity-introduction}
 Let $X$ be as above, and let $M$ be a closed surface, equipped with a Riemannian metric. If $u\colon M\to X$ is a continuous Sobolev map minimizing the energy in its homotopy class then $u$ is H\"older continuous. If $M$ is diffeomorphic to $S^2$ then every energy minimizer in a homotopy class is infinitesimally quasiconformal.
\et

In order to obtain the H\"older continuity in the result above we actually prove that energy minimizers in homotopy classes are harmonic, i.e.~locally energy minimizing. The H\"older regularity of harmonic maps in this context is known by \cite{LW16-harmonic}. The second part of the theorem holds more generally when $X$ is a complete metric space.  Combining \Cref{thm:full-Sacks-Uhlenbeck-metric-space} and \Cref{thm:regularity-introduction} we finally obtain:

\bc
 Let $X$ be as above. If $\pi_2(X)$ is non-trivial then $X$ contains a non-trivial infinitesimally quasiconformal harmonic sphere.
\ec

We end this subsection with some open questions. Using the techniques developed in \cite{SU-1981}, Sacks-Uhlenbeck moreover showed that every closed Riemannian manifold with non-trivial $k$-th homotopy group for some $k\geq 2$ receives a non-trivial harmonic $2$-sphere, see \cite[Theorem 5.7]{SU-1981}. We do not know whether a similar result holds in our metric setting and therefore ask:

\begin{quest}
    Let $X$ be a compact metric space with non-trivial $k$-th homotopy group for some $k\geq 2$. Under what additional conditions does $X$ admit a non-trivial harmonic $2$-sphere? 
\end{quest}

\Cref{thm:regularity-introduction} shows that energy minimizing spheres in homotopy classes are infinitesimally quasiconformal.  We do not know whether the same holds true in the more general context of harmonic spheres.

\begin{quest}
    Let $X$ be as in \Cref{thm:main-introduction} and let $u\colon S^2\to X$ be a harmonic map. Is it true that $u$ is infinitesimally quasiconformal?
\end{quest}

\subsection{Outline of Proofs}
We give insight into the proof of our main result, \Cref{thm:main-introduction}, by describing its two key ingredients. Let $X$ and $M$ be as in the statement of the theorem and let $\varphi\colon M\to X$ be continuous. The general idea is to apply a variational approach: consider an energy minimizing sequence $(u_n)$ of continuous Sobolev maps in the homotopy class of $\varphi$. By the Rellich-Kondrachov compactness theorem we may assume, after possibly passing to a subsequence, that $(u_n)$ converges in $L^2$ to a Sobolev map. This limit map need not have a continuous representative and even if it does, it might be in a different homotopy class than $\varphi$. Notice that this is unavoidable since, as described above, homotopy classes do not in general contain energy minimizers. Our first main ingredient, \Cref{thm:energy-minimizing-sequence}, provides a sufficient condition for the limit to be continuous and stay in the same homotopy class. We state a simplified version here.

\begin{ing}\label{ingredient-1}
    There exists $\varepsilon_0>0$ with the following property. Let $u$ be the $L^2$-limit of the energy minimizing sequence $(u_n)$. If there is $r_0>0$ such that for every $n$ the energy of $u_n$ is bounded by $\varepsilon_0$ on every $r_0$-ball then $u$ has  a continuous representative which is homotopic to $\varphi$.
\end{ing}

In particular, by the lower semicontinuity of the energy, the limit is an energy minimizer in the homotopy class of $\varphi$.

To prove \Cref{ingredient-1}, we use methods introduced in \cite{Soultanis-Wenger-2022} to show that we may perturb a fine triangulation of $M$ such that the following holds: the restriction of $u$ to the 1-skeleton $M^1$ is essentially continuous, and, up to taking a subsequence, $(u_n|_{M^1})$ has uniformly bounded length and converges uniformly to the continuous representative of $u|_{M^1}$. We want to compare $u$ to a continuous Sobolev map $v$ that agrees, on $M^1$, with the continuous representative of $u|_{M^1}$ and solves Dirichlet's problem on each 2-cell of the perturbed triangulation. Fix $n$ large enough. For every 2-cell $\Delta$, we use the local quadratic isoperimetric inequality of $X$ to construct a Sobolev annulus of small area connecting $u_n|_{\partial\Delta}$ and $v|_{\partial\Delta}$. As the sequence has small energy on $r_0$-balls, the Sobolev annulus glued to $u_n|_{\Delta}$ and $v|_{\Delta}$ results in a Sobolev sphere of small area, which is null-homotopic by \Cref{thm:small-area-spheres-contractible}. In particular, this induces a homotopy between $u_n$ and $v$. Lower semicontinuity of energy and the fact that $(u_n)$ is energy minimizing eventually imply that on each 2-cell, $u$ minimizes energy among Sobolev maps of the same trace. Therefore, $u$ has a continuous representative by \cite{LW16-harmonic}. As above, we can construct a homotopy between $u_n$ and the continuous representative of $u$. 

Our next goal is to find a condition on a continuous map $\varphi\colon M\to X$ ensuring the applicability of \Cref{ingredient-1} to an energy minimizing sequence in the homotopy class of $\varphi$.
We say that a map $\varphi\colon M\to X$ is $\varepsilon$-indecomposable for some $\varepsilon>0$ if all decompositions of $\varphi$ satisfy $e(\varphi_0)+e(\varphi_1)\geq e(\varphi)+\varepsilon$. Our second ingredient, found in \Cref{prop:unif-distributed-energy} and \Cref{prop:unif-distributed-energy-sphere}, implies the following.

\begin{ing}\label{ingredient-2}
    If $\varphi\colon M\to X$ is $\varepsilon$-indecomposable for some $0<\varepsilon<\varepsilon_0$, then there exists $r_0>0$ such that the following holds: if $u$ almost minimizes energy among Sobolev maps in the homotopy class of $\varphi$, then, up to precomposition with a conformal diffeomorphism, the map $u$ has energy bounded by $\varepsilon_0$ on every $r_0$-ball in $M$. 
\end{ing}

The precomposition with a conformal diffeomorphism is only needed when $M$ is diffeomorphic to $S^2$.

To prove \Cref{ingredient-2}, we let $r_0$ be well-chosen and assume that there is a ball $B\subset M$ of radius less than $r_0$ such that the energy of $u|_B$ is more than $\varepsilon_0$. We then find a curve $\gamma$ surrounding $B$ of small length and energy. As a consequence of the energy filling inequality proved in \cite[Section 4.1]{LW16-harmonic}, which depends on $X$ admitting a local quadratic isoperimetric inequality, we can fill the curve $\gamma$ with a continuous Sobolev map $v$ that has sufficiently small energy. We then use $v$ to construct a decomposition of $\varphi$ that contradicts the $\varepsilon$-indecomposability of $\varphi$. Note that the definition of decomposition is more restrictive in the case of $M$ being diffeomorphic to $S^2$. This case is dealt with by using a similar construction as above while precomposing $u$ with a suitable conformal diffeomorphism of the sphere to ensure that both maps in the decomposition of $\varphi$ are essential.

\subsection{Structure of article} 
In \Cref{sec:prelim}, we introduce relevant notation and concepts. In particular, we state the Reshetnyak definition of Sobolev spaces of metric space valued maps. The notion of general definitions of energy is recalled in \Cref{sec:energy-new}. The main result of \Cref{sec:cont-Sobolev-in-homotopy-class} shows that all homotopy classes of maps from a closed surface to $X$ contain a Sobolev map. In \Cref{sec:spheres-of-small-area}, we show that Sobolev spheres in $X$ of small area are null-homotopic. The purpose of \Cref{sec:regularity} is to establish that homotopic energy minimizers are harmonic and have Hölder continuous representatives. This together with the results from \Cref{sec:energy-new} imply \Cref{thm:regularity-introduction}. In \Cref{sec:1-homotopy-classes}, we recall relevant notions from \cite{Soultanis-Wenger-2022} and we prove an existence result of energy minimizers in given 1-homotopy classes, which provides an alternative proof of \Cref{cor:trivial-second-homotopy-group-intro}. In \Cref{sec:min-seq} and \Cref{sec:unif-distributed-energy}, we prove our two main ingredients described in the outline above. Finally, in \Cref{sec:proof-main-thm}, we prove our main result \Cref{thm:decomposition-minimizers} and its consequences.

\medskip

{\bf Acknowledgments:} Parts of this paper were finalized while the authors participated in the Trimester program {\it Metric Analysis} at the Hausdorff Research Institute for Mathematics in Bonn. We thank the institute for the hospitality and the inspiring atmosphere.

\section{Preliminaries}\label{sec:prelim}
\subsection{Notation}
Let $(X,d)$ be a metric space. The open ball in $X$ centered at $x\in X$ of radius $r>0$ will be denoted by $$B(x,r)=\{y\in X:d(x,y)<r\}.$$ The open unit disc in Euclidean $\R^2$ is referred to as $D$, the unit circle as $S^1$, and the standard 2-sphere in $\R^3$ as $S^2$. The length of a curve $\gamma$ in $X$ is denoted by $\ell(\gamma)$. We say that $X$ is called quasiconvex if there exists $\lambda\geq 1$ such that every pair of points $x,y\in X$ can be joined by a curve $\gamma$ satisfying $\ell(\gamma)\leq\lambda \cdot d(x,y)$. The space $X$ is geodesic if it is quasiconvex with constant $1$.

For $s\geq 0$, we denote the $s$-dimensional Hausdorff measure of a set $A\subset X$ by $\hm^s(A)$. The normalizing constant is chosen in such a way that if $X$ is the Euclidean space $\R^n$, then the Lebesgue measure agrees with $\hm^n$. If $M$ is a smooth manifold of dimension $n$ equipped with a Riemannian metric $g$, then the $n$-dimensional Hausdorff measure $\hm_g^n$ on $M$ coincides with the Riemannian volume. Recall from the introduction that, in this work, any closed, connected, smooth $2$-dimensional manifold $M$ is called a closed surface.

\subsection{Stereographic projections}\label{sec:stereographic-projection}
Throughout the paper we will need several basic facts about the stereographic projection based at the north pole.
Let $p_+$ and $p_-$ be the north and south pole of $S^2$, respectively.
The stereographic projection based at $p_+$ is the map $\psi\colon S^2\setminus\{p_+\}\to \C$ defined by $\psi(x,y,z) = \frac{1}{1-z}\cdot (x+i y)$. It satisfies 
\begin{equation}\label{eq:stereo-proj-ball-complement}
    \psi(S^2\setminus B(p_+, s)) = B(0, \cot(s/2))
\end{equation}
for all $s\in(0,\pi)$ and hence  $$\psi(B(p_-, t)) = B(0, \tan(t/2))$$  whenever $t\in(0,\pi)$. Denote by $\varrho\colon \C\to S^2\setminus \{p_+\}$ the inverse of $\psi$ and define $h(r) = 2\arctan(r)$. It follows from the above that 
\begin{equation}\label{eq:stereo-proj-inverse-ball-origin}
    \varrho(B(0,r))= B(p_-, h(r))
\end{equation}
for every $r>0$. Notice that we have $r\leq h(r)\leq 2r$ for all $0<r<1$. 

\subsection{Polyhedral complexes and triangulations} \label{subsec: Polyhedral_Complex_Triangulations}
A finite collection $K$ of compact convex polytopes, called cells of $K$, in some $\R^m$ is a polyhedral complex if each face of a cell is in $K$ and the intersection of two cells of $K$ is a face of each of them. The $j$-skeleton of $K$, denoted $K^j$, is the union of all cells of $K$ of dimension at most $j$. We also write $K$ to denote the union of all cells and  equip $K$ with the induced metric from $\R^m$, implying that an $n$-cell $\Delta$ of $K$ is isometric to a compact convex polytope in $\R^n$.

Let $M$ be a smooth closed $n$-manifold. A triangulation of $M$ consists of a polyhedral complex $K$ and a homeomorphism $h\colon K\to M$, where $h$ restricted to any $n$-cell $\Delta$ of $K$ is a $C^1$-diffeomorphism onto its image.

\subsection{Sobolev maps to metric spaces}\label{sec:Sobolev}
There exist several equivalent definitions of Sobolev maps from a Euclidean or Riemannian domain into a complete metric space, see for example \cite{Amb90}, \cite{KS93}, \cite{Res97}, \cite{Haj96}, \cite{HKST15}. We will review that of Reshetnyak \cite{Res97}. 

Let $(X, d)$ be a complete metric space and $M$ a smooth compact $m$-dimen\-sional Riemannian manifold, possibly with non-empty boundary. Let $g$ be the fixed Riemannian metric on $M$ and let $\Omega\subset M$ be an open set and $p>1$.

We denote by $L^p(\Omega, X)$ the collection of measurable and essentially separably valued maps $u\colon \Omega\to X$ such that for some and thus every $x\in X$ the function $u_x\colon\Omega\to \R$, defined by $$u_x(z)= d(x, u(z)),$$ belongs to the classical space $L^p(\Omega)$. A sequence $(u_n)\subset L^p(\Omega, X)$ is said to converge to $u\in L^p(\Omega, X)$ in $L^p(\Omega, X)$ if $$\int_\Omega d^p(u(z), u_n(z))\,d\hm^m_g(z)\to 0$$ as $n\to\infty$. 

\bd\label{def:Sobolev}
 A map $u\in L^p(\Omega, X)$ belongs to the Sobolev space $W^{1,p}(\Omega, X)$ if for every $x\in X$ the function $u_x$ belongs to $W^{1,p}(\Omega\setminus\partial M)$ and there exists $h\in L^p(\Omega)$ such that for all $x\in X$ we have $|\nabla u_x|_g\leq h$ almost everywhere on $\Omega$.  
\ed

In the definition above, $\nabla u_x$ is the weak gradient of $u_x$. Further, $|\cdot|_g$ denotes the norm induced by $g$. The Reshetnyak $p$-energy of $u\in W^{1,p}(\Omega, X)$ is defined by $$E_+^p(u)= \inf\left\{\|h\|_{L^p(\Omega)}^p: \text{$h$ as in \Cref{def:Sobolev}}\right\}.$$

Sobolev maps $u\in W^{1,p}(\Omega, X)$ have the following approximate metric differentiability property. For almost every $z\in \Omega$ there exists a unique seminorm $\apmd u_z$ on $T_zM$ such that $$\ap \lim_{v\to 0} \frac{d(u(\exp_z(v)), u(z)) - \apmd u_z(v)}{|v|_g}=0.$$ Here, $\exp_z\colon T_zM\to M$ denotes the exponential map, and $\ap\lim$ is the approximate limit, see \cite{Kir94} and e.g.\ \cite[Theorem 1.15 and Property 2.7]{Kar07}. If $\gamma$ is a Sobolev map defined on an interval or on $S^1$, we will write $|\gamma'|(t)$ instead of $\apmd \gamma_t(1)$.

As is shown in \cite[Page 1145]{LW15-Plateau}, the Reshetnyak energy satisfies $$E_+^p(u) = \int_\Omega\mathbf{I}_+^p(\apmd u_z)\,d\hm_g^m(z),$$ where for a seminorm $s$ on $\R^m$ we define $$\mathbf{I}_+^p(s) = \max\{s(v)^p: |v|=1\}$$ and extend it to seminorms on $T_zM$ by identifying $(T_zM, g(z))$ with $(\R^m, |\cdot|)$ via a linear isometry. For Sobolev maps $\gamma$ defined on an interval $(a, b)$ the Reshetnyak $p$-energy is denoted by $E^p$ and is simply given by $$E^p(\gamma)=\int_a^b|\gamma'|^p(t)\,dt,$$ and the length of the continuous representative of $\gamma$, denoted again by $\gamma$, satisfies $$\length(\gamma) = \int_a^b |\gamma'|(t)\,dt.$$

Next, we consider the case that $M$ has dimension 2. We define the Jacobian $\jac(s)$ of a seminorm $s$ on Euclidean $\R^2$ as the Hausdorff $2$-measure on $(\R^2,s)$ of the unit square if $s$ is a norm and zero otherwise. The definition extends to seminorms on $T_zM$ after identifying $(T_zM,g(z))$ with $(\R^2,|\cdot|)$ via a linear isometry. We can now define a notion of area of a Sobolev map.

\bd
 Let $u\in W^{1,2}(\Omega, X)$. The parametrized (Hausdorff) area of $u$ is defined by $$\Area(u)=\int_{\Omega} \jac(\apmd u_z)\,d\hm^2_g(z).$$
\ed

The parametrized area of a Sobolev map is invariant under precompositions with biLipschitz homeo\-morphisms, and thus independent of the Riemannian metric $g$.

Next we recall the concept of infinitesimal quasiconformality. This is a non-smooth substitute for the notion of weak conformality.

\bd\label{def:infinitesimal-qc}
    A map $u\in W^{1,2}(\Omega, X)$ is called infinitesimally $Q$-quasiconformal (with respect to a Riemannian metric $g$) if for almost every $z\in \Omega$ the seminorm $\apmd u_z$ satisfies
    $$\apmd u_z(v)\leq Q\cdot\apmd u_z(w)$$
    for all $v,w\in T_z M$ with $|v|_g=|w|_g$.
\ed

For every $u\in W^{1,2}(\Omega, X)$ it holds that $\jac(\apmd u_z)\leq \mathbf{I}_+^2(\apmd u_z)$, which we can integrate to obtain
\begin{align}\label{eq:area_less_than_energy}
{\rm Area}(u)\leq E^2_+(u).
\end{align}
If $u$ is infinitesimally $Q$-quasiconformal, then it follows from a short calculation that $\mathbf{I}_+^2(\apmd u_z)\leq Q^2\cdot \jac(\apmd \allowbreak u_z)$, and thus
$E^2_+(u)\leq Q^2\cdot{\rm Area}(u)$ in this case.

Finally, we recall the definition of the trace of a Sobolev map. Let $\Omega \subset M\setminus\partial M$ be a Lipschitz domain. Then for every $z$ in the boundary $\partial \Omega$ of $\Omega$ there exist an open neighborhood $U\subset M$ and a biLip\-schitz map $\psi\colon (0,1)\times [0,1)\to M$ such that $\psi((0,1)\times (0,1)) = U\cap \Omega$ and $\psi((0,1)\times\{0\}) = U\cap \partial\Omega$. Let $u\in W^{1,2}(\Omega, X)$. For almost every $s\in (0,1)$ the map $t\mapsto u\circ\psi(s,t)$ belongs to $W^{1,2}((0,1),X)$ by a Fubini-type argument and thus has an absolutely continuous representative which we denote by the same expression. The trace of $u$ is defined by $$\trace(u)(\psi(s,0))=\lim_{t\searrow 0} (u\circ\psi)(s,t)$$ for almost every $s\in(0,1)$. It can be shown that the trace is independent of the choice of the map $\psi$ and defines an element of $L^2(\partial \Omega, X)$, see \cite[Section 1.12]{KS93} for further details.

\section{General definitions of energy}\label{sec:energy-new}
In the setting of metric spaces there are several natural notions of energy and our results apply to general definitions of energy. The purpose of this short section is to recall the notion of a general definition of energy and to show that energy minimizers in homotopy classes are infinitesimally quasiconformal.

Let $\mathfrak{S}_2$ be the proper metric space of seminorms on $\R^2$ with the metric defined by $$d_{\mathfrak{S}_2}(s,s')=\max_{v\in S^1}|s(v)-s'(v)|.$$ The following definition appears in \cite{LW17-en-area}.

\bd
A definition of energy is a continuous map $\mathbf{I}\colon\mathfrak{S}_2\to[0,\infty)$ with the following properties.
\begin{enumerate}
    \item Monotonicity: $\mathbf{I}(s)\geq \mathbf{I}(s')$ whenever $s\geq s'$.
    \item Homogeneity: $\mathbf{I}(\lambda\cdot s)=\lambda^2\cdot \mathbf{I}(s)$ for all $\lambda\in[0,\infty)$.
    \item ${\rm SO}_2$-invariance: $\mathbf{I}(s\circ T)=\mathbf{I}(s)$ for any $T\in{\rm SO}_2$.
    \item Properness: The set $\mathbf{I}^{-1}([0,1])$ is compact in $\mathfrak{S}_2$.
\end{enumerate}
\ed

In addition, throughout this article, we always assume that a definition of energy $\mathbf{I}$ is quasiconvex, a condition that is defined below and ensures lower semicontinuity of energy. 

Some prominent definitions of energy are the Reshetnyak energy $\mathbf{I}^2_+$ as defined in the previous section and the Korevaar-Schoen energy $\mathbf{I}_{\text{avg}}$ given by
$$\mathbf{I}_{\text{avg}}(s)=\frac{1}{\pi}\int_{S^1}s(v)^2\,dv.$$
The latter agrees with the Dirichlet energy when $X$ is Riemannian, see \cite{KS93}. Properness and homogeneity imply that any two definitions of energy are comparable. In particular, for any given definition of energy $\mathbf{I}$, there exists $k_\mathbf{I}\geq1$ such that
\begin{equation}\label{eq:equivalence-energies}
    k_\mathbf{I}^{-1}\cdot \mathbf{I}\leq \mathbf{I}_+^2\leq k_\mathbf{I}\cdot \mathbf{I}.
\end{equation}

Let $M$ be a surface, equipped with a Riemannian metric $g$. Given an open set $\Omega\subset M$ the $\mathbf{I}$-energy of a Sobolev map $u\in W^{1,2}(\Omega,X)$ is defined by $$E_\mathbf{I}(u,g)=\int_{\Omega}\mathbf{I}(\apmd u_z)\,d\hm_g^2(z).$$ 

Note that $E_\mathbf{I}(u,g)$ is invariant under precompositions with conformal maps. If it is clear from the context which Riemannian metric $g$ we consider on $M$, we simply write $E_\mathbf{I}(u)$. 

The following theorem is a consequence of (the proofs of) \cite[Lemmas 3.1 and 4.1]{LW17-en-area}, compare with \cite[Theorem 6.2]{LW15-Plateau} and \cite[Theorem 6.6]{LW-param}.

\bt\label{thm:precomp-bilip-sphere}
 Let $X$ be a complete metric space and $\mathbf{I}$ a definition of energy. If $u\in W^{1,2}(S^2, X)$ satisfies $$E_{\mathbf{I}}(u)\leq E_{\mathbf{I}}(u\circ\eta)$$ for every biLipschitz homeomorphism $\eta$ of $S^2$ then $u$ is infinitesimally quasiconformal.
\et

The quasiconformality constant only depends on the definition of energy. When $\mathbf{I}$ is the Reshetnyak energy $\mathbf{I}_+^2$, then a stronger regularity property, called infinitesimal isotropy, holds for $u$, see \cite[Lemma 3.2]{LW17-en-area}. It implies for example that $u$ is infinitesimally $\sqrt{2}$-quasiconformal and when $X$ is a Riemannian manifold, or more generally a space with property (ET) (cf. \cite[Definition 11.1]{LW15-Plateau}), then infinitesimal isotropy is equivalent to weak conformality, thus $u$ is infinitesimally 1-quasi\-conformal, see \cite[Theorem 11.3]{LW15-Plateau}.

\bc\label{cor:energy-minimizer-inf-qc}
 Let $X$ be a complete metric space and $\mathbf{I}$ a definition of energy. If $u\in W^{1,2}(S^2, X)$ is continuous and minimizes energy in its homotopy class then $u$ is infinitesimally quasiconformal.
\ec

\begin{proof}
    Let $\eta \colon S^2\to S^2$ be a biLipschitz homeomorphism. After possibly precomposing it with an orientation reversing isometry of $S^2$ we may assume that $\eta$ is orientation preserving. In particular, $\eta$ is homotopic to the identity, see e.g.\ \cite[Page~51]{Milnor65}, and hence $u\circ\eta$ is homotopic to $u$, consequently $E_{\mathbf{I}}(u) \leq E_{\mathbf{I}}(u\circ\eta)$. Since $\eta$ was arbitrary it follows from \Cref{thm:precomp-bilip-sphere} that $u$ is infinitesimally quasiconformal.
\end{proof}

 The regularity properties of energy minimizers now allow us to show that existence of energy minimizers in a homotopy class is not guaranteed in the case when $\pi_2(X)$ is  non-trivial, even when $X$ satisfies strong assumptions.

\begin{figure}
    \centering
    \includegraphics[width=0.55\linewidth]{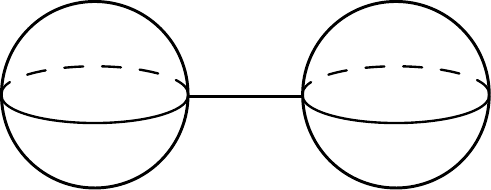}
    \caption{Illustration of the space $X$ in \Cref{example:dumbbell}}
    \label{fig:dumbbell}
\end{figure}

\begin{example}\label{example:dumbbell}
    Let $X$ be the metric space obtained by gluing two copies of $S^2$ to the endpoints of the interval $[0,1]$, and equipping it with the intrinsic length metric; for an illustration see \Cref{fig:dumbbell}. 
    Let $\varphi$ be a continuous map wrapping the sphere $S^2$ around both copies of $S^2$ in $X$. It is not hard to see that for any given energy $\mathbf{I}$, the homotopy class of $\varphi$ does not contain an energy minimizer. Indeed, assume $u\in W^{1,2}(S^2,X)$ is continuous and minimizes the $\mathbf{I}$-energy in the homotopy class of $\varphi$. Then it is infinitesimally quasiconformal by \Cref{cor:energy-minimizer-inf-qc}.  Since the approximate metric derivative of $u$ is degenerate for almost every $z\in U\coloneqq u^{-1}((0,1))$ by \cite[Proposition 4.3]{LW15-Plateau}, it then follows that $\apmd u_z$ is the zero-seminorm for almost every $z\in U$. Thus, the $\mathbf{I}$-energy of $u|_U$ must be zero, implying that $u|_U$ is locally constant, which is impossible.
\end{example}

Throughout this paper we will require that all our chosen definitions of energy $\mathbf{I}$ are quasiconvex in the following sense. For every finite-dimensional normed space $Y$ and every linear map $L\colon\R^2\to Y$ we have $$E_{\mathbf{I}}(L|_D)\leq E_\mathbf{I}(\nu)$$ for every smooth immersion $\nu\colon\bar{D}\to Y$ with $\nu|_{S^1}=L|_{S^1}$. The Reshetnyak definition of energy and the Korevaar-Schoen definition of energy are both quasiconvex, see the proofs of \cite[Corollaries 5.6 and 5.7]{LW15-Plateau}. Let $M$ be a closed surface, equipped with a Riemannian metric, and let $\Omega\subset M$ be an open subset. It follows from \cite[Theorem 5.4]{LW15-Plateau} that if $\mathbf{I}$ is quasiconvex then every sequence $(u_n)\subset W^{1,2}(\Omega, X)$ of uniformly bounded $\mathbf{I}$-energy and converging in $L^2(\Omega,X)$ to $u\in W^{1,2}(\Omega, X)$ satisfies $$E_\mathbf{I}(u)\leq\liminf_{n\to\infty}E_\mathbf{I}(u_n).$$

Lower semicontinuity of the $\mathbf{I}$-energy is crucial for proving existence of $\mathbf{I}$-energy minimizers. Such an existence result is provided by the following theorem, which will be repeatedly used throughout this work.

\bt\label{thm:existence-energy-min-given-trace}
 Let $X$ be a proper metric space admitting a local quadratic isoperimetric inequality and let $\Omega\subset M$ be a Lipschitz domain. Then for every Sobolev map $u\in W^{1,2}(\Omega, X)$ with bounded trace there exists $v\in W^{1,2}(\Omega, X)$ with $$E_{\mathbf{I}}(v) = \inf\{ E_{\mathbf{I}}(w): w\in W^{1,2}(\Omega, X), \ \trace(w) = \trace(u)\}$$ and such that $\trace(v) = \trace(u)$. Any such $v$ has a locally H\"older continuous representative $\bar{v}$. Moreover, if $\trace(u)$ is continuous then $\bar{v}$ extends continuously to the boundary.
\et

\begin{proof}
  The first part is a direct consequence of the Rellich-Kondrachov compactness theorem, see \cite[Theorem 1.13]{KS93}, together with  \cite[Lemma 4.11]{LW15-Plateau} and the lower semicontinuity of the $\mathbf{I}$-energy. The continuity (up to the boundary) follows from \cite[Theorem 1.3]{LW16-harmonic}.
\end{proof}

\section{Constructing continuous Sobolev maps in a homotopy class}\label{sec:cont-Sobolev-in-homotopy-class}
Let $\varphi\colon M\to X$ be a continuous map from a closed surface $M$ to a complete metric space $X$. We define the possibly empty family $$\Lambda(\varphi)\coloneqq \big\{u\in W^{1,2}(M,X): \text{$u$ continuous and homotopic to $\varphi$}\big\}.$$ The goal of the current section is to show that, in the setting of our main result, $\Lambda(\varphi)$ is non-empty.

\bp\label{prop:Sobolev-homotopic-to-cont}
 Let $X$ be a proper, quasiconvex metric space admitting a local quadratic isoperimetric inequality. Suppose furthermore that every continuous map from $S^2$ to $X$ of sufficiently small diameter is null-homotopic. Let $M$ be a closed surface. Then for every continuous map $\varphi\colon M\to X$ the family $\Lambda(\varphi)$ is non-empty.
\ep

The proposition will be derived from the following two lemmas.

\bl\label{lem:Sobolev-super-critical-approx-cont}
 Let $X$ be a proper, quasiconvex metric space admitting a local quadratic isoperimetric inequality. Then there exists $p>2$ such that for every closed surface $M$, every continuous map $\varphi\colon M\to X$, and every $\varepsilon>0$ there is $u\in W^{1,p}(M, X)$ continuous with $d(u(z), \varphi(z))\leq \varepsilon$ for every $z\in M$.
\el

\begin{proof}
  Let $\varphi\colon M\to X$ be a continuous map from a closed surface $M$ and let $\varepsilon>0$. Assume that $X$ satisfies a local quadratic isoperimetric inequality of constant $C>0$ and up to scale $l_0>0$. Let $\delta>0$ be sufficiently small, to be determined later, and fix a triangulation of $M$ such that each $2$-cell is a simplex. For every $j\in\{0,1,2\}$ we identify the $j$-skeleton of the triangulation with a subset of $M$ which we denote by $M^j$. Choosing the triangulation sufficiently fine, we may assume that $\diam(\varphi(\Delta))\leq \delta$ for every $2$-cell $\Delta\subset M^2$. For each vertex $a\in M^0$ set $u(a)= \varphi(a)$. If $e\subset M^1$ is an edge with endpoints $a$ and $b$ then let $u|_e$ be a constant speed Lipschitz curve of length at most $\lambda d(\varphi(a), \varphi(b))$, where $\lambda$ is the quasiconvexity constant of $X$. Fix a $2$-simplex $\Delta\subset M^2$ and let $\eta\colon \bar{D}\to \Delta$ be a biLipschitz homeomorphism such that the restriction to $S^1$ is a constant speed parametrization of $\partial\Delta$. Notice that $u\circ \eta|_{S^1}$ has length at most $3\lambda\delta$, which we may assume to be smaller than $ l_0$. By \cite[Theorem 3.4]{LW16-harmonic}, there exist $p>2$ and $C'$ only depending on $C$ and a Sobolev map $v\in W^{1,p}(D, X)$ with trace $u\circ\eta|_{S^1}$ and energy $$E_+^p(v)\leq C'E^p(u\circ\eta|_{S^1})\leq (2\pi)^{1-p}C'\cdot (3\lambda\delta)^p.$$ In particular, $v$ has a representative (denoted by $v$ again) that is $\alpha$-H\"older on all of $\bar{D}$ with $\alpha = 1-\frac{2}{p}$ and H\"older constant $C''\lambda\delta$ for some $C''$ only depending on $C$, see \cite[Proposition 3.3]{LW15-Plateau}. We have $\diam(v(\bar{D}))\leq 2C''\lambda\delta$. Now extend $u$ to $\Delta$ by setting $u=v\circ\eta^{-1}$ and notice that for each $z\in \Delta$ we have $$d(\varphi(z), u(z)) \leq \delta + 2C''\lambda\delta,$$ which we may assume to be smaller than $\varepsilon$ by choosing $\delta$ small enough. We may repeat this procedure for every $2$-cell $\Delta\subset M^2$, and thus, extend $u$ to all of $M$. Then $u$ is continuous and by the Sobolev gluing theorem, see \cite[Theorem 1.12.3]{KS93}, we conclude that $u\in W^{1,p}(M,X)$.
\end{proof}

We will furthermore need the following result, which is a consequence of the arguments in the proof of \cite[Proposition 6.2]{LWY20} together with \cite[Theorem 5.2 and Proposition 2.2]{LWY20}.

\bl\label{lem:curves-small-diam-null-homotopic}
 Let $X$ be a proper, quasiconvex metric space admitting a local quadratic isoperimetric inequality. Then there exist $r_0>0$ and $C_0\geq 1$ such that every closed curve with image in a ball $B(x,r)$ for some $r\in (0,r_0)$ is null-homotopic in $B(x, C_0r)$.
\el

\Cref{prop:Sobolev-homotopic-to-cont} is an almost direct consequence of the two lemmas above. We provide the proof for completeness. 

\begin{proof}[Proof of \Cref{prop:Sobolev-homotopic-to-cont}]
 Let $\delta_0>0$ be such that every continuous map from $S^2$ to $X$ of diameter at most $\delta_0$ is null-homotopic. Let $M$ be a closed surface and let $\varphi\colon M\to X$ be continuous. Let $\varepsilon>0$ be sufficiently small, to be determined later. By \Cref{lem:Sobolev-super-critical-approx-cont} there exists a continuous Sobolev map $u\in W^{1,2}(M, X)$ such that $d(\varphi(z), u(z))\leq \varepsilon$ for every $z\in M$. We now show that if $\varepsilon$ was chosen small enough, then $u$ is homotopic to $\varphi$ and thus $\Lambda(\varphi)$ is non-empty. 

 Indeed, fix a triangulation of $M$ such that each $2$-cell is a simplex. For every $j\in\{0,1,2\}$, we identify the $j$-skeleton of the triangulation with a subset of $M$ and write $M^j$. By choosing the triangulation fine enough, we may assume that $\varphi(\Delta)\cup u(\Delta)$ is contained in a ball of radius $2\varepsilon$ for every $2$-simplex $\Delta\subset M^2$. We build a homotopy $H\colon M\times[0,1]\to X$ from $\varphi$ to $u$ as follows. Firstly, we set $H(z,0) \coloneqq \varphi(z)$ and $H(z,1)\coloneqq u(z)$ for all $z\in M$. For each vertex $a\in M^0$, we let $H|_{\{a\}\times[0,1]}$ be a Lipschitz curve of length at most $\lambda\varepsilon$ from $\varphi(a)$ to $u(a)$, where $\lambda$ is the quasiconvexity constant of $X$. Thus, for each edge $e\subset M^1$, the curve $H|_{\partial(e\times[0,1])}$ is contained in a ball $B(x, 2(1+\lambda)\varepsilon)$. By choosing $\varepsilon>0$ so small that \Cref{lem:curves-small-diam-null-homotopic} is applicable, the map $H|_{\partial(e\times[0,1])}$ extends to a continuous map $H|_{e\times[0,1]}$ with image in $B(x, 2C_0(1+\lambda)\varepsilon)$, where $C_0$ is the constant from \Cref{lem:curves-small-diam-null-homotopic}. Then for each $2$-simplex $\Delta\subset M^2$, the map $H|_{\partial(\Delta\times[0,1])}$ has diameter at most $8C_0(1+\lambda)\varepsilon$, implying that $H|_{\partial(\Delta\times[0,1])}$ is null homotopic if $\varepsilon$ was chosen small enough such that $8C_0(1+\lambda)\varepsilon\leq\delta_0$. This means that we may extend $H$ continuously to its full domain $M\times[0,1]$ and thus $H$ is a homotopy from $\varphi$ to $u$.
\end{proof}

\section{Sobolev spheres of small area are null-homotopic}\label{sec:spheres-of-small-area}
In this section, we aim to prove that, under the same assumptions on $X$ as in the previous section, the condition of continuous spheres of small diameter in $X$ being null-homotopic implies that continuous Sobolev spheres in $X$ of small area are null-homotopic. Note that in general, a Sobolev sphere of small area can have large diameter.
\bt\label{thm:small-area-spheres-contractible}
 Let $X$ be a proper, quasiconvex metric space admitting a local quadratic isoperimetric inequality. Suppose furthermore that every continuous map from $S^2$ to $X$ of sufficiently small diameter is null-homotopic. Then there exists $\alpha_0>0$ such that every continuous map $u\in W^{1,2}(S^2, X)$ with $\Area(u)<\alpha_0$ is null-homotopic.
\et

For a metric space $X$ we denote by $\ell^\infty(X)$ the Banach space of all bounded real-valued functions equipped with the supremum norm. Recall that $X$ embeds isometrically into $\ell^\infty(X)$ by the Kuratowski embedding  $\iota\colon X\to\ell^{\infty}(X)$ at some fixed $x_0\in X$ defined by sending $x\in X$ to the function $$\iota(x)(\cdot)\coloneqq d(x,\cdot)-d(x_0,\cdot),$$ see \cite{Kuratowski1935}. The filling radius of a continuous map $\eta\colon S^2\to Y$ to a metric space $Y$ is defined by $$\Fillrad_Y(\eta)\coloneqq \inf\left\{R>0: \text{$\eta$ is null-homotopic in $N_R(\eta(S^2))$}\right\}.$$ Here, $N_R(A)$ denotes the open $R$-neighborhood of a set $A\subset Y$ for $R>0$. The proof of \Cref{thm:small-area-spheres-contractible} uses the following proposition.

\bp\label{prop:contracting-nbhd-linfty}
 Let $X$ be a complete metric space, and let $\iota\colon X\to Y\coloneqq\ell^\infty(X)$ be a Kuratowski embedding. Then for every continuous Sobolev map $u\in W^{1,2}(S^2, X)$ we have $$\Fillrad_Y(\iota\circ u)\leq L\cdot (\Area(u))^{\frac{1}{2}},$$ where $L$ is a universal constant.
\ep

\begin{proof}
 Set $v=\iota\circ u$ and let $\varepsilon>0$. By the proof of \cite[Proposition 3.1]{LWY20}, there exists a Lipschitz map $w\colon S^2\to Y$ with $\Area(w)\leq \Area(u) +\varepsilon$ and such that $$d(v(z), w(z))\leq \varepsilon$$ for all $z\in S^2$. It follows from the area formula for Lipschitz maps \cite[Corollary~8]{Kir94} that $$\hm^2(w(S^2))\leq \Area(w)\leq \Area(u)+\varepsilon.$$

 By \cite[Theorem 1.5]{Avvakumov-Nabutovsky}, there exists a Lipschitz map $H\colon S^2\times[0,1]\to Y$ such that $H_0=H(\cdot, 0)$ agrees with $w$, $H_1=H(\cdot, 1)$ has image in a $1$-dimensional simplicial complex $K\subset Y$, and $H(S^2\times[0,1])$ is contained in the $R$-neighborhood of $w(S^2)$, where $$R=L\cdot (\hm^2(w(S^2)))^{\frac{1}{2}}$$ for some universal constant $L$.

 As $H(S^2\times\{1\})$ is $1$-dimensional, the map $H_1$ is null-homotopic within its image by \cite[Corollary]{Curtis-Fort-1957}. From this and the above it follows that $v$ is null-homotopic in the $R'$-neighborhood of $v(S^2)$ in $Y$ for $$R'= \varepsilon+ L\cdot (\Area(u)+\varepsilon)^{\frac{1}{2}}.$$ Since $\varepsilon>0$ was arbitrary, the claim follows.
\end{proof}

With the aid of \Cref{prop:contracting-nbhd-linfty}, we may now prove the main result of this section.

\begin{proof}[Proof of \Cref{thm:small-area-spheres-contractible}]
Let $\delta_0>0$ be so small that every continuous map from $S^2$ to $X$ of diameter less than $\delta_0$ is null-homotopic. View $X$ as a subset of $Y\coloneqq \ell^\infty(X)$ via a Kuratowski embedding. Let $\alpha_0>0$ be small, to be determined later, and let $u\colon S^2\to X$ be a continuous Sobolev map with $\Area(u)<\alpha_0$. By \Cref{prop:contracting-nbhd-linfty}, there exists a homotopy $G\colon S^2\times[0,1]\to Y$ from $u$ to a constant map with $$G(S^2\times[0,1])\subset N_R(u(S^2)),$$ where $R = L\cdot(\Area(u))^{\frac{1}{2}} < L\sqrt{\alpha_0}$ for some universal constant $L$.

Fix a triangulation of $S^2\times[0,1]$ such that each cell is a simplex. For every $j\in\{0,1,2,3\}$ we identify the $j$-skeleton of the triangulation with a subset of $S^2\times[0,1]$, denoted $Z^j$. By choosing the triangulation sufficiently fine, we may assume that $$d(G(a),G(b))\leq \sqrt{\alpha_0}$$ whenever $a,b\in Z^0$ are endpoints of an edge $e\in Z^1$. Define a map $H\colon S^2\times [0,1]\to X$ as follows. Let $H$ agree with $G$ on $S^2\times\{0\}$, and define $H$ on $S^2\times\{1\}$ as a constant map with image equal to a nearest point of $G(S^2\times\{1\})$ in $u(S^2)$. For each vertex $a\in Z^0$ not contained in $S^2\times\{0\}$ or $S^2\times\{1\}$ let $H(a)$ be a nearest point of $G(a)$ in $u(S^2)$. In particular, if $a,b\in Z^0$ are endpoints of an edge $e\in Z^1$, then $$d(H(a), H(b))\leq (2L+1)\sqrt{\alpha_0}.$$ We may extend $H$ to $Z^1$ using the quasiconvexity of $X$ in a way that $$\diam(H(e))\leq (2L+1)\lambda\sqrt{\alpha_0}$$ for every edge $e\subset Z^1$, where $\lambda$ is the quasiconvexity constant of $X$. Assume that $\alpha_0$ was chosen small enough that \Cref{lem:curves-small-diam-null-homotopic} is applicable. Thus, for every $2$-simplex $\Delta\subset Z^2$ we can extend $H|_{\partial\Delta}$ to $\Delta$ such that $$\diam(H(\Delta))\leq 2(2L+1)\lambda C_0\sqrt{\alpha_0},$$ where $C_0\geq1$ is the constant from \Cref{lem:curves-small-diam-null-homotopic}. Finally, in case of $\alpha_0$ being so small that $4(2L+1)\lambda C_0\sqrt{\alpha_0}\leq\delta_0$, the map $H$ continuously extends to all of $S^2\times[0,1]$, which finishes the proof.
 \end{proof}

The following consequence of \Cref{thm:small-area-spheres-contractible} will be useful later.

 \bc\label{cor:homotopic-relative-boundary}
  Let $X$ and $\alpha_0$ be as in \Cref{thm:small-area-spheres-contractible}, and let $u, v\colon \bar{D}\to X$ be continuous and such that $u|_{S^1} = v|_{S^1}$. If $u$ and $v$ are Sobolev and satisfy $$\Area(u) + \Area(v) <\alpha_0$$ then $u$ and $v$ are homotopic relative to $S^1$.
 \ec

 \begin{proof}
     Let $\psi_-\colon S^2_-\to \bar D$ be the restriction of the stereographic projection to the lower half sphere $S_-^2$, and define $\psi_+\colon S^2_+\to \bar D$ by $\psi_+=\psi_-\circ  \sigma$, where $S_+^2$ is the upper hemisphere and $\sigma$ is given by $\sigma(x,y,z) = (x,y,-z)$. Notice that $\psi_-$ and $\psi_+$ are biLipschitz and agree on the common boundary. Therefore, the map $w\colon S^2\to X$ coinciding with $u\circ\psi_+$ on $S^2_+$ and with $v\circ\psi_-$ on $S^2_-$ is continuous. Moreover, by the Sobolev gluing theorem, see \cite[Theorem 1.12.3]{KS93}, $w$ is in $W^{1,2}(S^2, X)$ and satisfies
     $$
     \Area(w) = \Area(u\circ\psi_+)+\Area(u\circ\psi_-)=\Area(u) + \Area(v) < \alpha_0,
     $$
     since the area of a Sobolev map is invariant under precompositions with biLip\-schitz maps. It therefore follows from \Cref{thm:small-area-spheres-contractible} that there exists a continuous extension $\bar w$ 
    of $w$ to the closed $3$-dimensional unit ball. Let $\rho=\psi_-^{-1}$ and define $H\colon\bar D\times [0,1]\to X$ by $$H(z,t) = \bar w(t\rho(z)+(1-t)\sigma(\rho(z))).$$
     It is straight-forward to check that $H$ is a homotopy from $u$ to $v$ relative to $S^1$.
 \end{proof}

\section{Regularity of homotopic energy minimizing maps} \label{sec:regularity}
The aim of this section is to determine the main regularity properties of energy minimizers in homotopy classes. In addition, an example is given to illustrate that these properties are generally best possible.

Let $\mathbf{I}$ be a definition of energy as in \Cref{sec:energy-new}.

\bd
Let $X$ be a complete metric space, and let $M$ be a closed surface, equipped with a Riemannian metric. A map $u\in W^{1,2}(M,X)$ is called $\mathbf I$-harmonic if every $z\in M$ has an open neighborhood $U$ such that for each Lipschitz domain $\Omega$ with compact closure in $U$ we have $$E_{\mathbf I} (u|_\Omega)\leq E_{\mathbf{I}}(v)$$ for all $v\in W^{1,2}(\Omega,X)$ with $\trace(v) = \trace(u|_\Omega)$.
\ed

The following proposition shows how $\mathbf{I}$-harmonicity is related to H\"older continuity.

\bp\label{prop:regularity}
 Let $M$ be a closed surface and $X$ a proper metric space admitting a local quadratic isoperimetric inequality. If $u\in W^{1,2}(M, X)$ is $\mathbf{I}$-harmonic then $u$ has a H\"older continuous representative.
\ep

The H\"older exponent only depends on $\mathbf{I}$ and the isoperimetric constant of $X$.

\begin{proof}
The definition of $\mathbf{I}$-harmonic along with inequality \eqref{eq:equivalence-energies} imply that $u$ is locally quasi-harmonic in the sense of \cite[Definition~1.1]{LW16-harmonic}. It follows from \cite[Theorem~1.3]{LW16-harmonic} that $u$ has a locally $\alpha$-Hölder continuous representative $\bar{u}$ for some $\alpha$ only depending on $\mathbf{I}$ and the isoperimetric constant of $X$. By compactness of $M$, the map $\bar{u}$ is globally $\alpha$-Hölder continuous.
\end{proof}

Together with \Cref{cor:energy-minimizer-inf-qc} the following theorem directly implies \Cref{thm:regularity-introduction}, the regularity result stated in the introduction.

\bt\label{thm:regularity}
 Let $M$ be a closed surface and $X$ a proper, quasiconvex metric space admitting a local quadratic isoperimetric inequality. Suppose furthermore that every continuous map from $S^2$ to $X$ of sufficiently small diameter is null-homotopic. If $u\in W^{1,2}(M, X)$ is continuous and minimizes energy in its homotopy class then $u$ is $\mathbf{I}$-harmonic and, in particular, has a H\"older continuous representative.
\et

\begin{proof}
    Let $\alpha_0>0$ be as in \Cref{thm:small-area-spheres-contractible}, and let $k_\mathbf{I}\geq1$ be the constant from \eqref{eq:equivalence-energies}. Fix a point $z\in M$ and let $ U\subset M$ be a disc with smooth boundary such that $z\in U$ and 
    \begin{equation}\label{ineq:energy-bound-regularity}
        E_\mathbf{I}(u|_{U})<\frac{\alpha_0}{2k_\mathbf{I}}. 
    \end{equation} 
    Let $v\in W^{1,2}(U,X)$ be an $\mathbf{I}$-energy minimizer in the class of all Sobolev maps in $W^{1,2}(U,X)$ having trace $\trace(u|_{U})$ as in \Cref{thm:existence-energy-min-given-trace}. In particular, we may assume that $v$ is continuous and extends continuously to the boundary. By the $\mathbf{I}$-energy minimality of $v$, and after applying \eqref{eq:area_less_than_energy},  \eqref{eq:equivalence-energies} and \eqref{ineq:energy-bound-regularity}, we obtain that $$\max\{\Area(u|_{U}),\Area(v)\}< \alpha_0/2.$$ Consequently, the mappings $u|_{\bar{U}}$ and $v$ are homotopic relative to $\partial U$ by \Cref{cor:homotopic-relative-boundary}. In particular, the map $w$ agreeing with $v$ on $U$ and with $u$ on $M\setminus {U}$ is homotopic to $u$. Note that by the Sobolev gluing theorem, see \cite[Theorem 1.12.3]{KS93}, we also have $w\in W^{1,2}(M,X)$. As $u$ minimizes the $\mathbf{I}$-energy in its homotopy class, it has to hold that $E_\mathbf{I}(u)\leq E_\mathbf{I}(w)$. In particular, by definition of $w$ and the $\mathbf{I}$-energy minimality of $v$, we obtain $E_\mathbf{I}(u|_U)=E_\mathbf{I}(v)$. Hence, $u$ minimizes $\mathbf{I}$-energy on $U$ and thus also on every Lipschitz domain contained in $U$. This shows that $u$ is $\mathbf{I}$-harmonic.
\end{proof}

Hölder regularity in \Cref{prop:regularity} and \Cref{thm:regularity} is the best we can hope for as illustrated by the following example.

\begin{example}
    Let $X$ be the geodesic metric space (biLipschitz homeomorphic to $S^2$) obtained by gluing two copies of a $2$-dimensional cone of small cone angle along their boundaries. Let $u\colon S^2\to X$ be the homeomorphism that maps each hemisphere to one of the cones via a radial map with radial stretch function $t\mapsto t^\alpha$ for some $\alpha\in(0,1)$. Then $u$ is $\alpha$-H\"older (and no better), and it can be shown that if $\alpha$ is chosen suitably (depending on the cone angle) then $u$ is an energy minimizer in its homotopy class, compare with \cite[Example 8.3]{LW15-Plateau}. 
\end{example}

\section{Minimizers in 1-homotopy classes and a proof of \texorpdfstring{\Cref{cor:trivial-second-homotopy-group-intro}}{Corollary \ref{cor:trivial-second-homotopy-group-intro}}}\label{sec:1-homotopy-classes}
In this section we review some of the constructions from \cite{Soultanis-Wenger-2022} that will be needed in the subsequent sections. In particular, we will recall the definitions of admissible deformation of a triangulation on a closed surface and the notion of $1$-homotopy type of a Sobolev map. We then show the existence of energy minimizers in a given $1$-homotopy class and deduce from this \Cref{cor:trivial-second-homotopy-group-intro} from the introduction. 

Let $M$ be a closed surface. An admissible deformation on $M$ is a smooth map ${\Phi\colon M\times\R^m\to M}$, $m\in\N$, where $\Phi_{\xi}\coloneqq\Phi(\cdot,\xi)$ is a diffeomorphism for every $\xi\in\R^m$ and $\Phi_0={\rm id}_M$, and such that the derivative of $\Phi^z\coloneqq\Phi(z,\cdot)$ satisfies $$D\Phi^z(0)(\R^m)=T_zM$$ for every $z\in M$. Admissible deformations exist on any compact surface, see \cite[Proposition 2.2]{Soultanis-Wenger-2022} and also \cite{White86}, \cite{White88}, \cite{HL03} for related results.

Let $h\colon K\to M$ be a triangulation of $M$. Fix an admissible deformation $\Phi\colon M\times\R^m\to M$ and, for $\xi\in\R^m$, denote by $h_{\xi}\colon K\to M$ the triangulation given by $h_{\xi}\coloneqq\Phi_{\xi}\circ h$. Let $X$ be a complete metric space and $u\in W^{1,2}(M, X)$. By \cite[Proposition 3.5]{Soultanis-Wenger-2022} there exists a ball $B_{\Phi,h}\subset\R^m$, centered at the origin and with radius only depending on $\Phi$ and $h$, such that $u\circ h_{\xi}|_{K^1}$ is essentially continuous for almost every $\xi\in B_{\Phi, h}$. Moreover, if $X$ is proper and admits a local quadratic isoperimetric inequality then, by \cite[Theorem 3.6]{Soultanis-Wenger-2022}, there is a negligible set $N\subset B_{\Phi, h}$ such that for all $\xi, \zeta\in B_{\Phi, h}\setminus N$ the continuous representatives of $u\circ h_{\xi}|_{K^1}$ and $u\circ h_{\zeta}|_{K^1}$, denoted by the same symbols,  are homotopic as maps from $K^1$ to $X$. The map $u$ is called $1$-homotopic to a given continuous map $\varphi\colon M\to X$ if $u\circ h_{\xi}|_{K^1}$ is homotopic to $\varphi\circ h|_{K^1}$ for almost every $\xi\in B_{\Phi, h}$. Moreover, two maps $u, v\in W^{1,2}(M, X)$ are called $1$-homotopic if $u\circ h_\xi|_{K^1}$ and $v\circ h_{\xi}|_{K^1}$ are homotopic for almost every $\xi\in B_{\Phi, h}$. If $X$ is, in addition, quasiconvex then being $1$-homotopic is independent of the choice of triangulation $h$ and of admissible deformation $\Phi$, see \cite[Theorem 4.1]{Soultanis-Wenger-2022}.

Let $\mathbf{I}$ be a definition of energy as in \Cref{sec:energy-new}. Using results from \cite{Soultanis-Wenger-2022} we can now prove the following theorem.

\bt\label{thm:1-homotopy}
 Let $X$ be a compact, quasiconvex metric space admitting a local quadratic isoperimetric inequality, and let $M$ be a closed surface, equipped with a Riemannian metric. Then for every continuous map $\varphi\colon M\to X$ there exists an $\mathbf{I}$-energy minimizer $u$ in the $1$-homotopy class of $\varphi$; each such $u$ is $\mathbf{I}$-harmonic and thus has a H\"older continuous representative.
\et

We emphasize that (the continuous representative) of $u$ need not be in the homotopy class of $\varphi$. However, $u$ is an energy minimizer in its own homotopy class since homotopic maps are in particular $1$-homotopic.

\begin{proof}
Let $\varphi\colon M\to X$ be continuous and denote by $\Lambda_1(\varphi)$ the family of Sobolev maps $1$-homotopic to $\varphi$. We first show that $\Lambda_1(\varphi)$ is non-empty. For this, let $\varepsilon>0$ be sufficiently small. By \Cref{lem:Sobolev-super-critical-approx-cont}, there exists a continuous Sobolev map $v\in W^{1,2}(M, X)$ satisfying $d(v(z), \varphi(z))\leq \varepsilon$ for all $z\in M$. Let $h\colon K\to M$ be a triangulation. It follows from the quasiconvexity of $X$ and from \Cref{lem:curves-small-diam-null-homotopic} that if $\varepsilon$ was chosen small enough, then $v\circ h|_{K^1}$ is homotopic to $\varphi\circ h|_{K^1}$. 
This easily implies that $v$ is $1$-homotopic to $\varphi$ and hence $v\in \Lambda_1(\varphi)$.

Next, let $(u_n)\subset \Lambda_1(\varphi)$ be an $\mathbf{I}$-energy minimizing sequence. After possibly passing to a subsequence we may assume, by the Rellich-Kondrachov compactness theorem, that $(u_n)$ converges in $L^2(M, X)$ to a Sobolev map $u\in W^{1,2}(M, X)$. It follows from \cite[Theorem 4.7]{Soultanis-Wenger-2022} that $u_n$ is $1$-homotopic to $u$ for every sufficiently large $n$ and hence $u$ is $1$-homotopic to $\varphi$. Consequently, $u\in\Lambda_1(\varphi)$ and by lower semicontinuity of the $\mathbf{I}$-energy, $u$ is an $\mathbf{I}$-energy minimizer in $\Lambda_1(\varphi)$.

Finally, let $u\in\Lambda_1(\varphi)$ be any $\mathbf{I}$-energy minimizer and let $z\in M$. Choose a triangulation $h\colon K\to M$ such that $z$ is contained in $U=h(\Delta)$ for some open cell $\Delta\subset K$ and let $\Omega\subset U$ be a Lipschitz domain compactly contained in $U$. Then $\Omega\subset h_\xi(\Delta)$ for all sufficiently small $\xi\in B_{\Phi, h}$. In particular, if $v\in W^{1,2}(\Omega, X)$ satisfies $\trace(v) = \trace(u|_\Omega)$ then the Sobolev map $w\in W^{1,2}(M, X)$ which agrees with $v$ on $\Omega$ and with $u$ on $M\setminus \Omega$ is $1$-homotopic to $u$. Consequently, $E_{\mathbf{I}}(w)\geq E_{\mathbf{I}}(u)$ and thus $E_{\mathbf{I}}(v)\geq E_{\mathbf{I}}(u|_\Omega)$, which shows that $u$ is $\mathbf{I}$-harmonic. The existence of a H\"older continuous representative of $u$ now follows from \Cref{prop:regularity}. 
\end{proof}

When $X$ has trivial second homotopy group then continuous $1$-homotopic maps are homotopic, see \cite[Lemma 6.2]{Soultanis-Wenger-2022}, and thus we conclude the following statement that implies \Cref{cor:trivial-second-homotopy-group-intro}.
 
\bc
Let $X$ and $M$ be as in \Cref{thm:1-homotopy}. If $\pi_2(X)$ is trivial, then every continuous map from $M$ to $X$ contains an $\mathbf{I}$-energy minimizer in its homotopy class, and any such energy minimizer is $\mathbf{I}$-harmonic.
\ec

We moreover obtain the following result which generalizes \cite[Theorem 2.1]{SY79} and \cite[Theorem 5.2]{SU-1981}. Note that each free homotopy class of (unbased) continuous maps from $M$ to $X$ induces a conjugacy class of homomorphisms from $\pi_1(M)$ to $\pi_1(X)$.
\bc
 Let $X$ and $M$ be as in \Cref{thm:1-homotopy}. Then every conjugacy class of homomorphisms between $\pi_1(M)$ and $\pi_1(X)$ is induced by a continuous Sobolev map from $M$ to $X$ that minimizes the $\mathbf{I}$-energy among all continuous Sobolev maps inducing the same conjugacy class.
\ec

\begin{proof}
Let $f\colon\pi_1(M)\to \pi_1(X)$ be a representative of a conjugacy class of homomorphisms. Note that since $M$ is a surface, there exists a based continuous map $\varphi\colon M\to X$ such that it induces $f$ as a homomorphism between the fundamental groups. Let $\Lambda$ denote the space of continuous Sobolev maps inducing a conjugacy class containing $f$. By the first part of the proof of \Cref{thm:1-homotopy} there exists a continuous Sobolev map $1$-homotopic to $\varphi$ and this map induces the same conjugacy class of homomorphisms as $\varphi$. Hence $\Lambda$ is non-empty. 

Let $e\coloneqq\inf\{E_\mathbf{I}(u): u\in \Lambda\}$ and let $(u_n)$ be an $\mathbf I$-energy minimizing sequence in $\Lambda$. Then by passing to a subsequence we may assume, by the Rellich-Kondrachov compactness theorem, that $(u_n)$ converges in $L^2(M,X)$ to a Sobolev map $u\in W^{1,2}(M,X)$. By lower semicontinuity of the $\mathbf{I}$-energy we have $E_\mathbf{I}(u)\leq e$. It follows from \cite[Theorem 4.7]{Soultanis-Wenger-2022} that $u_n$ is 1-homotopic to $u$ for some sufficiently large $n$. By \Cref{thm:1-homotopy} there exists an $\mathbf{I}$-energy minimizer $v$ that is $1$-homotopic to $u_n$, and each such is Hölder continuous. It follows that $v\in \Lambda$, and since $e\leq E_\mathbf{I}(v)\leq E_\mathbf{I}(u) \leq e$, we conclude that $v$ is an $\mathbf{I}$-energy minimizer in $\Lambda$.
\end{proof}

\section{Minimizing sequences of uniformly distributed energy}\label{sec:min-seq}
The purpose of this section is to establish a convergence result for energy minimizing sequences that have uniformly distributed energy.

Let $M$ be a closed surface equipped with a Riemannian metric. Let $X$ be a compact, quasiconvex metric space admitting a local quadratic isoperimetric inequality. Suppose furthermore that every continuous map from $S^2$ to $X$ of sufficiently small diameter is null-homotopic. Let $\mathbf{I}$ be a definition of energy as in \Cref{sec:energy-new}. Given a continuous map $\varphi\colon M\to X$, we define $$e_\mathbf{I}(\varphi)\coloneqq\inf\left\{E^2_\mathbf{I}(u):u\in\Lambda(\varphi)\right\}.$$ Set $\varepsilon_\mathbf{I}\coloneqq k_\mathbf{I}^{-1}\alpha_0$, where $k_\mathbf{I}$ is the constant from \eqref{eq:equivalence-energies} and $\alpha_0>0$ is as in \Cref{thm:small-area-spheres-contractible}. Under these assumptions we will establish the following theorem, which is one of the crucial ingredients in the proof of the main result of our work.

\bt\label{thm:energy-minimizing-sequence}
Let $(u_n)$ be a sequence of continuous maps in $W^{1,2}(M,X)$ converging in $L^2(M,X)$ to a map $u\in W^{1,2}(M,X)$. Suppose that
\begin{equation}\label{eq:energy-minimizing-sequence}
    E_\mathbf{I}(u_n)-e_\mathbf{I}(u_n)\to 0\quad\text{as }n\to\infty,
\end{equation}
and that there exists $r_0>0$ such that
\begin{equation}\label{eq:local-energy-bound}
    E_\mathbf{I}(u_n|_{B(p,r_0)})\leq\frac{\varepsilon_\mathbf{I}}{5}
\end{equation}
for every $p\in M$ and all $n\in\N$. Then $u$ has a continuous representative $\bar u$ satisfying $E_\mathbf{I}(\bar{u})=e_\mathbf{I}(\bar{u})$, and $u_n$ is homotopic to $\bar u$ for sufficiently large $n$.
\et

In the proof of this theorem we will need the two lemmas below.

\bl\label{lma:sobolev-homotopy-of-small-area}
Let $X$ be a proper, quasiconvex metric space admitting a local quadratic isoperimetric inequality of constant $C$ and up to scale $l_0$. Denote by $\lambda$ the quasiconvexity constant of $X$ and let $0<\varepsilon<\frac{ l_0}{4\lambda}$. Assume that $\gamma_0,\gamma_1\in W^{1,2}([0,1], X)$ are continuous and such that $$d(\gamma_0(t), \gamma_1(t))\leq \varepsilon$$ for every $t\in [0,1]$. Let $\beta_0,\beta_1\in W^{1,2}([0,1], X)$ be continuous with $\ell(\beta_j)<\lambda\varepsilon$, and assume that $\beta_j$ connects $\gamma_0(j)$ and $\gamma_1(j)$ for $j\in\{0,1\}$. Then there exists a continuous map $\chi\in W^{1,2}([0,1]^2,X)$ with $\chi(\cdot,j)=\gamma_j$ and $\chi(j,\cdot)=\beta_j$ for $j\in\{0,1\}$ and such that $$\Area(\chi)\leq C'\varepsilon\cdot \max\{\varepsilon, \length(\gamma_0), \length(\gamma_1)\},$$ where $C'$ only depends on $C$ and $\lambda$.
\el

\begin{proof}
Set $\ell_1\coloneqq \max\{\varepsilon, \length(\gamma_0), \length(\gamma_1)\}$ and let $k\in\N$ be the smallest integer such that $\ell_1/k\leq \varepsilon$. Choose partitions $0=t_0^j<t_1^j<\dots<t_k^j=1$ of $[0,1]$ satisfying $$\length(\gamma_j|_{[t_m^j, t_{m+1}^j]}) = \frac{\length(\gamma_j)}{k}$$ for $j=0,1$ and $m=0,1,\dots, k-1$. We denote by $0=t_0<t_1<\dots<t_n=1$ the joint partition, in particular, $n\leq 2k$.

We set $\tilde{\beta}_0=\beta_0$ and $\tilde{\beta}_n=\beta_1$ and, for $i\in\{1,\dots, n-1\}$, let $\tilde{\beta}_i\colon[0,1]\to X$ be a Lipschitz curve from $\gamma_0(t_i)$ to $\gamma_1(t_i)$ of length $\ell(\tilde{\beta}_i)\leq \lambda\cdot d(\gamma_0(t_i),\gamma_1(t_i))\leq\lambda\cdot\varepsilon$. The map $\chi_i\colon \partial([t_i,t_{i+1}]\times[0,1])\to X$ defined by setting $\chi_i(\cdot,j)=\gamma_j$ for $j\in\{0,1\}$ and $\chi_i(j,\cdot)=\tilde{\beta}_j$ for $j\in\{t_i,t_{i+1}\}$ is continuous and piecewise Sobolev with length $$\length(\chi_i)\leq \frac{\length(\gamma_0)}{k} + \frac{\length(\gamma_1)}{k} + 2\lambda\varepsilon\leq \frac{2\ell_1}{k} + 2\lambda\varepsilon\leq 4\lambda\varepsilon< l_0.$$ Therefore, by the local quadratic isoperimetric inequality, \cite[Lemma 2.6]{LW16-harmonic} and \cite[Proposition 2.4]{Soultanis-Wenger-2022}, the map $\chi_i$ extends to a continuous map $\Bar{\chi}_i\colon [t_i, t_{i+1}]\times[0,1]\to X$ which is Sobolev and satisfies $$\Area(\bar\chi_i)\leq 2C\cdot\ell(\chi_i)^2\leq 32 C\lambda^2\varepsilon^2.$$ Gluing the $\Bar{\chi}_i$ gives a continuous map $\chi\in W^{1,2}([0,1]^2,X)$, by the Sobolev gluing theorem, with $\chi(\cdot,j)=\gamma_j$ and $\chi(j,\cdot)=\beta_j$ for $j\in\{0,1\}$ and which satisfies $$\Area(\chi)\leq n\cdot 32C\lambda^2\varepsilon^2\leq 64C\lambda^2k\varepsilon^2\leq 128C\lambda^2\cdot \ell_1\cdot \varepsilon.$$ This concludes the proof.
\end{proof}

In the next lemma, we use the notation introduced in \Cref{sec:1-homotopy-classes}. Let $M$ be a closed surface, equipped with a Riemannian metric. Fix a triangulation $h\colon K\to M$ and let $\Phi\colon M\times \R^m\to M$ an admissible deformation on $M$. 

\bl \label{lem:technical-lemma-1-skeleton-regularity}
Let $X$ be a complete metric space and let $(u_n)\subset W^{1,2}(M, X)$ be a sequence with uniformly bounded energy converging in $L^2(M, X)$ to some $u\in W^{1,2}(M, X)$. Then for almost every $\xi\in B_{\Phi, h}$ the following holds:
\begin{enumerate}
   \item the restriction $u\circ h_\xi|_{K^1}$ is essentially continuous and $u\circ h_\xi|_e$ is in $W^{1,2}(e, X)$ for every edge $e\subset K^1$; moreover $$\trace(u\circ h_\xi|_{\Delta}) = u\circ h_\xi|_{\partial \Delta}$$ for every $2$-simplex $\Delta\subset K^2$.
   \item there exists a subsequence $(u_{n_j})$ such that the continuous representative $\gamma_j$ of $u_{n_j}\circ h_\xi|_{K^1}$ converges uniformly to the continuous representative of $u\circ h_\xi|_{K^1}$; moreover, for each edge $e\subset K^1$ the lengths $\ell(\gamma_j|_e)$ are uniformly bounded.
\end{enumerate}
\el

The lemma can be proved by combining arguments from the proofs of \cite[Lemmas 2.6 and 3.7]{Soultanis-Wenger-2022}. For the convenience of the reader we provide a detailed proof. We will use the notion of Newton-Sobolev maps. For the definition we refer for example to \cite{Soultanis-Wenger-2022}.

\begin{proof}
By \cite[Lemma 3.3]{Soultanis-Wenger-2022} there exists $L>0$ such that for every Borel function $\rho\colon M\to[0,\infty]$ we have 
\begin{equation}\label{eq:integral-admissible-deformation}
  \int_{B_{\Phi, h}}\int_{K^1}\rho^2\circ h_\xi(z)\,d\hm^1(z)\,d\xi\leq L \int_M\rho^2\,d\hm^2.
\end{equation}
Let $u\in W^{1,2}(M,X)$ and let $v\colon M\to X$ be a Newton-Sobolev representative of $u$ with upper gradient $\rho\in L^2(M)$. In particular, for every Lipschitz curve $\beta\colon [a,b]\to M$ we have $$d(v(\beta(a)), v(\beta(b)))\leq \int_a^b \rho(\beta(s))|\beta'(s)|\,ds,$$ see \cite[Proposition 2.5]{Soultanis-Wenger-2022}. By \cite[Corollary 3.4]{Soultanis-Wenger-2022}, for almost every $\xi$ the maps $v\circ h_\xi$ and $u\circ h_\xi$ agree almost everywhere on $K^1$. This shows that it suffices to prove (i) with $u$ replaced by $v$. Since $\rho\in L^2(M)$ it follows from \eqref{eq:integral-admissible-deformation} that
\begin{equation}\label{eq:int-on-1-skeleton-finite}
\int_{K^1}\rho^2\circ h_\xi(z)\,d\hm^1(z)<\infty   
\end{equation}
for almost every $\xi\in B_{\Phi, h}$. Fix such $\xi$, let $e\subset K^1$ be a closed edge, and let $x,y\in e$ be distinct points. If $\gamma\colon [a,b]\to e$ is a Lipschitz curve connecting $x$ and $y$ then 
\begin{equation*}
  \begin{split}
      d(v\circ h_\xi(x), v\circ h_\xi(y))&\leq \int_a^b\rho(h_\xi\circ\gamma(s)) |(h_\xi\circ\gamma)'(s)|\,ds\\
      &\leq \lip(h_\xi)\int_a^b\rho\circ h_\xi(\gamma(s))|\gamma'(s)|\,ds.
  \end{split}
\end{equation*}
This together with \eqref{eq:int-on-1-skeleton-finite} implies that $v\circ h_\xi|_e$ is Newton-Sobolev and in particular in $W^{1,2}(e, X)$. If $\gamma$ is an arc-length parametrization of the segment from $x$ to $y$ then the above together with H\"older's inequality implies $$d(v\circ h_\xi(x), v\circ h_\xi(y))\leq   \sqrt{A}\lip(h_\xi)\cdot |x-y|^{\frac{1}{2}},$$ where $A= \int_e\rho^2\circ h_\xi(z)\,d\hm^1(z)$. In particular, $v\circ h_\xi|_e$ is H\"older continuous. One shows analogously that 
\begin{equation*}
 \ell(v\circ h_\xi|_e)\leq \sqrt{A}\lip(h_\xi)\cdot \ell(e).
\end{equation*}
Finally, it follows from \cite[Proposition 2.5]{Soultanis-Wenger-2022} that for almost every $\xi\in B_{\Phi, h}$ and for every $2$-simplex $\Delta\subset K^2$ we have $\trace(u\circ h_\xi|_\Delta) = v\circ h_\xi|_{\partial \Delta}$. This shows (i).

We turn to (ii) and let $(u_n)\subset W^{1,2}(M, X)$ be a sequence with uniformly bounded energy converging in $L^2(M, X)$ to $u$. By passing to a subsequence we may asssume that $u_n$ converges almost everywhere to $u$. For each $n$, let $v_n\colon M\to X$ be a Newton-Sobolev representative of $u_n$ and $\rho_n$ an upper gradient of $v_n$ with $$\|\rho_n\|^2_{L^2(M)} \leq 2 E_+^2(u_n).$$ By the first part of the proof and by \cite[Corollary 3.4]{Soultanis-Wenger-2022} there exists a negligible set $N\subset B_{\Phi, h}$ such that for every $\xi\in B_{\Phi, h}\setminus N$ and each $n\in\N$ the map $v_n\circ h_\xi|_{K^1}$ is the continuous representative of $u_n\circ h_\xi|_{K^1}$ and $v_n\circ h_\xi(z) \to u\circ h_\xi(z)$ for $\hm^1$-almost every $z\in K^1$. Combine \eqref{eq:integral-admissible-deformation} with Fatou's lemma to obtain
\begin{equation*}
  \begin{split}
      \int_{B_{\Phi, h}}\Big(\liminf_{n\to\infty}\int_{K^1}\rho_n^2\circ h_\xi(z)\,d\hm^1(z)\Big)\,d\xi&\leq \liminf_{n\to\infty} \int_{B_{\Phi, h}}\int_{K^1}\rho_n^2\circ h_\xi(z)\,d\hm^1(z)\,d\xi\\
      &\leq L\cdot \liminf_{n\to\infty} \int_M\rho_n^2\,d\hm^2 <\infty.
  \end{split}
\end{equation*}
In particular, for almost every $\xi\in B_{\Phi, h}\setminus N$ there exists a subsequence satisfying $$\sup_{j\in\N}\Big(\int_{K^1}\rho_{n_j}^2\circ h_\xi(z)\,d\hm^1(z)\Big)<\infty.$$ It follows from this together with  the first part of the proof that the sequence $(v_{n_j}\circ h_\xi|_e)$ is uniformly H\"older and has uniformly bounded length for each edge $e\subset K^1$. Consequently, $v_{n_j}\circ h_\xi|_{K^1}$ converges uniformly to the continuous representative of $u\circ h_\xi|_{K^1}$, and this completes the proof of (ii). 
\end{proof}

\begin{proof}[Proof of \Cref{thm:energy-minimizing-sequence}]
Let $\Phi\colon M\times\R^m\to M$ be an admissible deformation, and let $h\colon K\to M$ be a triangulation of $M$ where the 2-cells of $K$ are simplices and such that $\diam(h(\Delta))<\frac{r_0}{2}$ for every 2-cell $\Delta\subset K^2$. Recall the notation from \Cref{sec:1-homotopy-classes}. 

Let $(u_n)$ and $u$ be as in the statement of the theorem and notice that \eqref{eq:local-energy-bound} implies that the energy of $u_n$ is uniformly bounded. By \Cref{lem:technical-lemma-1-skeleton-regularity} there exist $\xi\in B_{\Phi,h}$ and a subsequence $(u_{n_j})$ with the following properties. Firstly, $u\circ h_\xi|_{K^1}$ is essentially continuous. We denote by $u\circ h_\xi|_{K^1}$ again the continuous representative. Secondly, $u\circ h_\xi|_e$ is Sobolev for every edge $e\subset K^1$ and $\trace(u\circ h_\xi|_\Delta) = u\circ h_\xi|_{\partial \Delta}$ for every $2$-simplex $\Delta\subset K^2$. Finally, for every edge $e\subset K^1$ the curves $u_{n_j}\circ h_\xi|_e$ are Sobolev, have uniformly bounded length, and converge uniformly to $u\circ h_\xi|_e$. We may of course assume that $\xi$ is so small that $\diam(h_\xi(\Delta))<r_0$ for every $2$-cell $\Delta\subset K^2$.

By \Cref{thm:existence-energy-min-given-trace} there exists, for each simplex $\Delta\subset K^2$, an $\mathbf{I}$-energy minimizer on $h_\xi(\Delta)$ with trace $u|_{h_\xi(\partial \Delta)}$ which can moreover be chosen to be continuous up to the boundary. By gluing all these maps, we obtain a continuous map $v\colon M\to X$ that is moreover Sobolev by the Sobolev gluing theorem. Notice that $$E_{\mathbf{I}}(v|_{h_\xi(\Delta)}) \leq E_{\mathbf{I}}(u|_{h_\xi(\Delta)}) \leq \liminf_{j\to\infty} E_{\mathbf{I}}(u_{n_j}|_{h_\xi(\Delta)}) \leq \frac{\varepsilon_I}{5}$$ for every simplex $\Delta\subset K^2$. 

We claim that $u_{n_j}$ is homotopic to $v$ for every sufficiently large $j$. Before proving the claim we first explain how it can be used to finish the proof. Since $v$ is homotopic to $u_{n_j}$, it follows from the above together with \eqref{eq:energy-minimizing-sequence} that $$e_{\mathbf{I}}(v) \leq E_{\mathbf{I}}(v)\leq E_{\mathbf{I}}(u)\leq \liminf_{j\to\infty} E_{\mathbf{I}}(u_{n_j}) = e_{\mathbf{I}}(v),$$ and hence equality holds everywhere. In particular, $E_{\mathbf{I}}(u|_{h_\xi(\Delta)}) = E_{\mathbf{I}}(v|_{h_\xi(\Delta)})$ for every simplex $\Delta\subset K^2$. Hence, $u|_{h_\xi(\Delta)}$ is an energy minimizer and thus has a representative which is continuous up to the boundary by \Cref{thm:existence-energy-min-given-trace}. Furthermore, $u|_{h_\xi(\Delta))}$ is homotopic to $v|_{h_\xi(\Delta)}$ relative to the boundary by \Cref{cor:homotopic-relative-boundary}. Therefore, the continuous representative $\bar{u}$ of $u$ is homotopic to $v$, and hence to $u_{n_j}$ for all $j$ large enough, and $E_{\mathbf{I}}(\bar{u}) = e_{\mathbf{I}}(\bar{u})$. It can easily be seen that actually $\bar{u}$ must be homotopic to $u_n$ for all $n$ large enough. This proves the theorem assuming the claim.

We are left to prove the claim. For this, we will construct a homotopy $H\colon K\times [0,1]\to X$ from $u_{n_j}\circ h_\xi$ to $v\circ h_\xi$ for all sufficiently large $j$. Let $l_0,\lambda,C'$ be as in \Cref{lma:sobolev-homotopy-of-small-area}, and fix $0<\varepsilon<\frac{ l_0}{4\lambda}$ small enough such that $$C'\varepsilon\cdot\max\left\{\varepsilon,\ell(v\circ h_\xi|_e), \ell(u_{n_j}\circ h_\xi|_e)\right\} \leq \frac{\alpha_0}{15}$$ for all edges $e\subset K^1$ and for all $j\in\N$. Then fix $j\in \N$ large enough so that $$d\left(v\circ h_\xi(x), u_{n_j}\circ h_\xi(x)\right)\leq\varepsilon$$ for all $x\in K^1$. Firstly, define $H(x,0)\coloneqq u_{n_j}\circ h_\xi(x)$ and $H(x,1) \coloneqq v\circ h_\xi(x)$ for all $x\in K$. Then, for each vertex $a\in K^0$, use the quasiconvexity of $X$ to define $H$ on $\{a\}\times[0,1]$ such that $H|_{\{a\}\times[0,1]}$ is a Lipschitz curve of length less than $\lambda\varepsilon$ from $u_{n_j}\circ h_\xi(a)$ to $ v\circ h_\xi(a)$. Since $u_{n_j}\circ h_\xi|_e$ and $v\circ h_\xi|_e$ are Sobolev curves for each edge $e\subset K^1$, we may apply \Cref{lma:sobolev-homotopy-of-small-area} to extend $H$ to the interior of $e\times[0,1]$ such that $H|_{e\times[0,1]}$ is a continuous Sobolev map satisfying $$\Area(H|_{e\times[0,1]})\leq C'\varepsilon\cdot\max\left\{\varepsilon,\ell(v\circ h_\xi|_e), \ell(u_{n_j}\circ h_\xi|_e)\right\}\leq \frac{\alpha_0}{15}.$$ Further, since our triangulation was chosen fine enough, we may apply \eqref{eq:local-energy-bound}, together with \eqref{eq:area_less_than_energy} and \eqref{eq:equivalence-energies}, to find that for all 2-cells $\Delta\subset K^2$
$$
\Area(u_{n_j}\circ h_\xi|_\Delta) = \Area(u_{n_j} |_{h_\xi(\Delta)}) \leq k_\mathbf{I}E_\mathbf{I}(u_{n_j}|_{h_\xi(\Delta)})\leq \frac{k_I\varepsilon_{\mathbf{I}}}{5}\leq \frac{\alpha_0}{5},
$$
and that $\Area(v\circ h_\xi|_\Delta)\leq \alpha_0/5$. Note that for each 2-cell $\Delta\subset K^2$ there exists a biLipschitz map $\mu_\Delta\colon S^2\to \partial (\Delta\times [0,1])$. By the Sobolev gluing theorem, it follows that $H\circ\mu_\Delta\in W^{1,2}(S^2,X)$. Note that every 2-cell $\Delta\subset K^2$ is a simplex, in particular $\partial\Delta$ has three edges $e_1,e_2, e_3$, and we can compute that
$$
\Area(H\circ\mu_\Delta) = \Area(u_{n_j}\circ h_\xi|_\Delta) + \Area( v\circ h_\xi|_\Delta) + \sum_{i=1}^3 \Area(H|_{e_i\times[0,1]}) \leq \frac{3\alpha_0}{5}. 
$$
Hence by \Cref{thm:small-area-spheres-contractible}, for each 2-cell $\Delta\subset K^2$, we can find a continuous extension of $H\circ \mu_\Delta$ to the closed unit ball in $\R^3$, that we then can use to define $H$ on the interior of $\Delta\times [0,1]$ in a way that $H|_{\Delta\times [0,1]}$ is continuous. After this last step, $H$ is a continuous map defined on all of $K\times [0,1]$, and it follows that $u_{n_j}\circ h_\xi$ is homotopic to $v\circ h_\xi$, which in particular implies that $u_{n_j}$ is homotopic to $v$. This proves the claim and completes the proof of the theorem.
\end{proof}

\section{Uniformly distributing energy}\label{sec:unif-distributed-energy}
The aim of this section is to find a condition ensuring that \Cref{thm:energy-minimizing-sequence} is applicable. Let $\mathbf{I}$ be a definition of energy as in \Cref{sec:energy-new} and let $M$, $X$ and $\varepsilon_\mathbf{I}$ be as in \Cref{sec:min-seq}. 

Given $\varepsilon>0$, we say that a continuous map $\varphi\colon M\to X$ is $\varepsilon$-indecomposable if $$e_\mathbf{I}(\varphi_0) + e_\mathbf{I}(\varphi_1) \geq e_\mathbf{I}(\varphi) + \varepsilon$$ for every decomposition $\varphi_0$ and $\varphi_1$ of $\varphi$. 

\bp\label{prop:unif-distributed-energy}
Suppose that $M$ is not diffeomorphic to a sphere. If  $\varphi\colon M\to X$ is $\varepsilon$-indecomposable for some $0<\varepsilon<\varepsilon_\mathbf{I}$, then there exists $r_0>0$ such that for every $u\in \Lambda(\varphi)$ with
$E_\mathbf{I}(u)\leq e_\mathbf{I}(\varphi)+10^{-1}\varepsilon$ we have
$$E_\mathbf{I}(u|_{B(p,r_0)})\leq\frac{\varepsilon_\mathbf{I}}{5}$$ 
for every $p\in M$.
\ep

In the proof of the proposition we will use the following consequence of the energy filling inequality proved in \cite[Section 4.1]{LW16-harmonic} together with \eqref{eq:equivalence-energies}. Let $l_0>0$ be the scale up to which the local isoperimetric inequality holds in $X$. Then every continuous curve $\gamma\in W^{1,2}(S^1,X)$ with $\ell(\gamma)<l_0$ is the trace of a map $w\in W^{1,2}(D, X)$ satisfying 
\begin{equation}\label{eq:energy-filling-inequality-I-energy}
    E_{\mathbf{I}}(w)\leq C_{\mathbf{I}}\cdot E^2(\gamma),
\end{equation}
 where $C_{\mathbf{I}}$ is a constant only depending on $\mathbf{I}$ and on the isoperimetric constant of $X$.

\begin{proof}[Proof of \Cref{prop:unif-distributed-energy}]
Let $\kappa\geq1$ be such that every $p\in M$ admits a conformal $\kappa$-biLipschitz chart $\varrho\colon D\to M$ with $\varrho(0)=p$. Let $0<\varepsilon<\varepsilon_{\mathbf{I}}$ and let $\varphi\colon M\to X$ be $\varepsilon$-indecomposable. We set 
 \begin{align*}
     \delta\coloneqq\min\left\{\frac{ l_0^2}{2\pi},\frac{\varepsilon}{10\,C_{\mathbf{I}}}\right\},\quad L\coloneqq e^{\delta^{-1}k_\mathbf{I}(e_\mathbf{I}(\varphi)+10^{-1}\varepsilon)},\quad\text{and}\quad r_0 \coloneqq \frac{1}{\kappa L},
 \end{align*}
 where $C_{\mathbf{I}}$ is as in \eqref{eq:energy-filling-inequality-I-energy} and $k_\mathbf{I}$ is as in \eqref{eq:equivalence-energies}. 
 
Let $u\in \Lambda(\varphi)$ satisfy $E_\mathbf{I}(u)\leq e_\mathbf{I}(\varphi)+10^{-1}\varepsilon$ and assume by contradiction that there exists $p\in M$ such that $$E_\mathbf{I}(u|_{B(p,r_0)})>\frac{\varepsilon_\mathbf{I}}{5}.$$ We first find a curve around $p$ along which $u$ has small energy. To achieve this, let $\varrho$ be a conformal $\kappa$-biLipschitz chart around $p$ as at the beginning of the proof. Since $B(p,r_0)\subset \varrho(B(0, \kappa r_0))$ the Sobolev map $v\coloneqq u\circ\varrho$ satisfies $$E_{\mathbf{I}}(v|_{B(0, \kappa r_0)}) = E_{\mathbf{I}}(u|_{\varrho(B(0, \kappa r_0))}) \geq E_{\mathbf{I}}(u|_{B(p, r_0)}) >\frac{\varepsilon_{\mathbf{I}}}{5}.$$ For $0<r<1$ define $\gamma_r\colon S^1\to D$ by $\gamma_r(z)= r\cdot z$. Then for almost every $r\in(0,1)$ the map $v\circ\gamma_r$ is Sobolev and $$|(v\circ \gamma_r)'|(t) = \apmd v_{\gamma_r(t)} (\gamma'_r(t))$$ for almost every $t\in S^1$, hence $$|(v\circ\gamma_r)'|^2(t) \leq r^2\cdot \mathbf{I}_+^2(\apmd v_{\gamma_r(t)})\leq r^2  k_{\mathbf{I}} \cdot\mathbf{I}(\apmd v_{\gamma_r(t)}).$$ Integrating in polar coordinates yields $$E_{\mathbf{I}}(v|_{D\setminus B(0,\kappa r_0)}) = \int_{\kappa r_0}^1\int_{S^1}r \cdot\mathbf{I}(\apmd v_{\gamma_r(t)})\,dt\,dr\geq k_{\mathbf{I}}^{-1}\cdot \int_{\kappa r_0}^1 r^{-1} E^2(v\circ\gamma_r)\,dr.$$ Using this together with the definition of $r_0$, we obtain the existence of a measurable subset $A\subset (\kappa r_0,1)$ of positive measure such that for every $r\in A$ the curve $v\circ\gamma_r$ is Sobolev and satisfies 
\begin{gather*}
    E^2(v\circ\gamma_{r}) \leq \frac{k_\mathbf{I}E_{\mathbf{I}}(v|_{D\setminus B(0,\kappa r_0)})}{\log(L)}<  \frac{k_\mathbf{I}(e_\mathbf{I}(\varphi)+10^{-1}\varepsilon)}{\log(L)}=\delta.
\end{gather*}
Fix $r\in A$ and notice that, by H\"older's inequality, we have $$\ell(v\circ\gamma_r)\leq \sqrt{2\pi \cdot E^2(v\circ\gamma_r)}<\sqrt{2\pi\delta} \leq l_0.$$ Let $\delta_r\colon \bar{D}\to \bar{B}(0,r)$ be the scaling map given by $\delta_r(z) = r\cdot z$ and notice that $\trace(v\circ\delta_r)=v\circ\gamma_r$. By \Cref{thm:existence-energy-min-given-trace} there exists $w\in W^{1,2}(D, X)$ with trace $v\circ\gamma_r$ and which minimizes the $\mathbf{I}$-energy among Sobolev maps with the same trace; moreover, $w$ is continuous and continuously extends to the boundary. The energy filling inequality \eqref{eq:energy-filling-inequality-I-energy} implies that $$E_{\mathbf{I}}(w) \leq C_{\mathbf{I}} \cdot E^2(v\circ\gamma_r) < C_{\mathbf{I}}\cdot\delta \leq \frac{\varepsilon}{10}.$$

We now use $w$ to define a decomposition of $u$ as follows. Let $B\subset M$ be the disc given by $B=\varrho(B(0,r))$. The map $u_0\colon M\to X$ agreeing with $u$ on $M\setminus B$ and with $\hat{w}\coloneqq w\circ (\varrho\circ \delta_r)^{-1}$ on $B$ is continuous and Sobolev by the Sobolev gluing theorem. Moreover, the inequalities above and the conformal invariance of the $\mathbf{I}$-energy yield $$E_{\mathbf{I}}(u_0)=  E_{\mathbf{I}}(u) - E_{\mathbf{I}}(u|_B) + E_{\mathbf{I}}(w) < E_{\mathbf{I}}(u) - \frac{\varepsilon}{10}\leq e_{\mathbf{I}}(\varphi)$$ and hence $u_0$ is not homotopic to $\varphi$. In particular, $u_0$ is not homotopic to $u$ and thus $\hat{w}$ is not homotopic to $u|_B$ relative to $\partial B$. Next, choose conformal diffeomorphisms $\eta_\pm\colon S^2_\pm\to \bar{B}$ agreeing on the common boundary and let $u_1\colon S^2\to X$ be the continuous map such that $u_1 = u\circ\eta_+$ on $S^2_+$ and $u_1=\hat{w}\circ\eta_-$ on $S^2_-$. Here, $S_+^2$ and $S_-^2$ denote the upper and lower hemisphere, respectively. By the above, $u_1$ is essential and thus $u_0$ and $u_1$ form a decomposition of $u$. Since $u$ is homotopic to $\varphi$ this decomposition is easily seen to give rise to a decomposition $\varphi_0$ and $\varphi_1$ of $\varphi$ such that $\varphi_0$ is homotopic to $u_0$ and $\varphi_1$ is homotopic to $u_1$. As $u_1$ is Sobolev and $$E_{\mathbf{I}}(u_1) = E_{\mathbf{I}}(u|_B) + E_{\mathbf{I}}(w)$$ we obtain $$e_{\mathbf{I}}(\varphi_0)+ e_{\mathbf{I}}(\varphi_1) \leq E_{\mathbf{I}}(u_0) + E_{\mathbf{I}}(u_1)= E_{\mathbf{I}}(u) + 2E_{\mathbf{I}}(w) < e_{\mathbf{I}}(\varphi) + \frac{3\varepsilon}{10},$$ which contradicts the $\varepsilon$-indecomposability of $\varphi$.
\end{proof}

We next prove a suitable version of the proposition above in the case that $M$ is homeomorphic to the sphere $S^2$. Since every Riemannian metric on $S^2$ is conformally equivalent to the standard metric, we may and will assume in the following that $S^2$ is equipped with the standard metric.

\bp\label{prop:unif-distributed-energy-sphere}
    If $\varphi\colon S^2\to X$ is $\varepsilon$-indecomposable for some $0<\varepsilon<\varepsilon_\mathbf{I}$, then there exists $r_0>0$ with the following property. For every $u\in \Lambda(\varphi)$ with $E_\mathbf{I}(u)\leq e_\mathbf{I}(\varphi)+10^{-1}\varepsilon$
    there is a conformal diffeomorphism $\eta\colon S^2\to S^2$ such that $$E_\mathbf{I}(u\circ\eta|_{B(p, r_0)})\leq \frac{\varepsilon_\mathbf{I}}{5}$$ for every $p\in S^2$.
\ep

The general strategy of proof is similar to that of \Cref{prop:unif-distributed-energy}. In addition we need the following simple fact.

\bl\label{lem:good-conf-diffeo-sphere}
 Let $L\geq 1$ and $0<r<L^{-1}$ and $p_0\in S^2$. Then there exist $0<r_0<L^{-1}$ only depending on $L$ and a conformal diffeomorphism $\eta\colon S^2\to S^2$ such that for every $p\in S^2$ we have $$\eta(B(p, r_0))\subset B(\eta(p), r) \quad \text{or} \quad \eta(B(p, r_0))\subset S^2\setminus B(p_0, Lr).$$
\el

 As the proof will show, the value $r_0$ decreases as the value of $L$ increases.

\begin{proof}
Let $L$, $r$, and $p_0$ be as in the statement of the lemma. We may assume that $p_0 = p_-$ is the south pole. Let $\psi$ be the stereographic projection based at the north pole $p_+$ and let $\varrho$ be its inverse.
Set $r_1\coloneqq \arccot(L)$ and let $\kappa$ denote the biLipschitz constant of the restriction of $\psi$ to $S^2\setminus B(p_+, r_1)$. Finally, set $r_0\coloneqq 2^{-1}\cdot\min\{r_1, L^{-1}, \kappa^{-2}\}$. 

Let $\delta_r\colon\C\to\C$ be the scaling map given by $\delta_r(z)=r\cdot z$ and denote by $\eta\colon S^2\to S^2$ the conformal diffeomorphism agreeing with $\varrho\circ\delta_r\circ\psi$ on $S^2\setminus\{p_+\}$. It follows from \eqref{eq:stereo-proj-ball-complement} and \eqref{eq:stereo-proj-inverse-ball-origin} that $$\eta(B(p_+, 2r_1)) = S^2\setminus B(p_-, h(Lr))\subset S^2\setminus B(p_-, Lr)$$ since $h(Lr)\geq Lr$. 

Now, let $p\in S^2$. We distinguish two cases. If $B(p, r_0)$ intersects $B(p_+, r_1)$ non-trivially then $B(p, r_0)$ is contained in $B(p_+, 2r_1)$ and hence, by the above, we have $\eta(B(p, r_0))\subset S^2\setminus B(p_-, Lr)$. If $B(p, r_0)$ does not intersect $B(p_+, r_1)$ then we have
\begin{equation*}
    \eta(B(p, r_0)) 
    \subset B(\eta(p), r\kappa^2r_0)
    \subset B(\eta(p), r)
\end{equation*}
since the restriction of $\psi$ to $S^2\setminus B(p_+, r_1)$ is $\kappa$-biLipschitz.
\end{proof}

\begin{proof}[Proof of \Cref{prop:unif-distributed-energy-sphere}]
Let $0<\varepsilon<\varepsilon_{\mathbf{I}}$ and let $\varphi\colon S^2\to X$ be an  $\varepsilon$-indecom\-posable map. Define, as in the proof of \Cref{prop:unif-distributed-energy}, $$\delta\coloneqq \min\left\{\frac{ l_0^2}{2\pi}, \frac{\varepsilon}{10 C_{\mathbf{I}}}\right\}  \quad\text{and}\quad L\coloneqq 2e^{\delta^{-1}k_\mathbf{I}(e_\mathbf{I}(\varphi) + 10^{-1}\varepsilon)},$$ where $l_0>0$ is the scale up to which the isoperimetric inequality holds in $X$, $C_{\mathbf{I}}$ is as in \eqref{eq:energy-filling-inequality-I-energy}, and $k_\mathbf{I}$ is as in \eqref{eq:equivalence-energies}. Let $0<r_0<L^{-1}$ be as in \Cref{lem:good-conf-diffeo-sphere}.

Let $u\in\Lambda(\varphi)$ be as in the statement of the proposition. For each $p\in S^2$ set $$r(p) \coloneqq \inf\left\{r>0: E_\mathbf{I}(u|_{B(p, r)})\geq 5^{-1}\varepsilon_\mathbf{I}\right\}$$ and let $\bar{r}>0$ be the infimum of the $r(p)$ over all $p\in S^2$. We clearly have $$E_\mathbf{I}(u|_{B(p,\bar{r})})\leq \frac{\varepsilon_\mathbf{I}}{5}$$ for every $p\in S^2$ and there exists $\bar{p}\in S^2$ such that equality holds for $p=\bar{p}$. If $\bar{r}\geq r_0$ then the proposition holds with $\eta$ being the identity mapping, so we may assume that $\bar{r}<r_0$. We claim that
\begin{equation}\label{eq:energy-bigger-ball-bounded-below}
    E_\mathbf{I}(u|_{B(\bar{p}, L\bar{r})})> E_\mathbf{I}(u) - \frac{\varepsilon_\mathbf{I}}{5}.
\end{equation}
The proposition easily follows from this together with \Cref{lem:good-conf-diffeo-sphere}. Indeed, let $\eta\colon S^2\to S^2$ be as in the lemma applied with $p_0=\bar{p}$ and $r=\bar{r}$. Then for every $p\in S^2$ we have $$\eta(B(p, r_0))\subset B(\eta(p), \bar{r})\quad \text{or} \quad \eta(B(p, r_0))\cap B(\bar{p}, L\bar{r})=\emptyset.$$ In the first case we obtain $$E_\mathbf{I}(u\circ\eta|_{B(p, r_0)})\leq E_\mathbf{I}(u|_{B(\eta(p), \bar{r})}) \leq \frac{\varepsilon_\mathbf{I}}{5}$$ and in the second case $$E_\mathbf{I}(u\circ\eta|_{B(p, r_0)})\leq E_\mathbf{I}(u|_{S^2\setminus B(\bar{p}, L\bar{r})}) \leq \frac{\varepsilon_\mathbf{I}}{5}.$$ This establishes the proposition assuming \eqref{eq:energy-bigger-ball-bounded-below}.

We are left to prove inequality \eqref{eq:energy-bigger-ball-bounded-below}. We argue by contradiction and assume that this inequality is false. Let $\varrho\colon D\to S^2$ be the conformal biLipschitz chart obtained by composing the restriction of the inverse stereographic projection to $D$ with an isometry of $S^2$ in such a way that $\varrho(0)= \bar{p}$. Notice that $$B(\bar{p}, r)\subset \varrho(B(0,r)) \subset B(\bar{p}, 2r)$$ for every $0<r<1$, see \Cref{sec:stereographic-projection}. Set $v\coloneqq u\circ\varrho$ and observe that $E_{\mathbf{I}}(v|_{B(0, \bar{r})})\geq \varepsilon_{\mathbf{I}}/5$ and $E_{\mathbf{I}}(v|_{B(0,2^{-1}L\bar{r})})\leq E_{\mathbf{I}}(u) - \varepsilon_{\mathbf{I}}/5$. Let $\gamma_r\colon S^1\to D$ be the curve given by $\gamma_r(z)= r\cdot z$. By comparing to the proof of \Cref{prop:unif-distributed-energy}, we can find the following estimate:
$$E_\mathbf{I}(v|_{B(0, 2^{-1}L \bar{r})\setminus B(0, \bar{r})})\geq k_{\mathbf{I}}^{-1}\int_{\bar{r}}^{2^{-1}L \bar{r}}r^{-1} E^2(v\circ\gamma_r)\,dr.$$ This implies the existence of a measurable subset $A\subset (\bar{r}, 2^{-1}L\bar{r})$ of positive measure such that for every $r\in A$ the curve $v\circ\gamma_r$ is Sobolev with $$E^2(v\circ\gamma_r)< \frac{k_{\mathbf{I}}(e_{\mathbf{I}}(\varphi)+10^{-1}\varepsilon)}{\log(L/2)}= \delta;$$ in particular, we have $\ell(v\circ\gamma_r)< \sqrt{2\pi\delta}\leq l_0$. Fix $r\in A$. Now, we have found a good estimate on the energy of the curve $v\circ\gamma_r$, in that it allows us to use the construction from the proof of \Cref{prop:unif-distributed-energy}. More specifically, we construct, exactly as in the proof of \Cref{prop:unif-distributed-energy}, the maps $u_0, u_1\colon S^2\to X$ that satisfy $E_\mathbf{I}(u_0)+E_\mathbf{I}(u_1)<e_\mathbf{I}(\varphi)+\varepsilon$, and such that $E_\mathbf{I}(u_0)<e_\mathbf{I}(\varphi)$. Further, our assumption that inequality  \eqref{eq:energy-bigger-ball-bounded-below} is false implies that $E_\mathbf{I}(u_1)<e_\mathbf{I}(\varphi)$.

In particular, $u_0$ and $u_1$ cannot be homotopic to $u$ and from this we deduce that they must both be essential and therefore form a decomposition of $u$. Consequently, they induce a decomposition $\varphi_0$ and $\varphi_1$ of $\varphi$ such that $\varphi_0$ is homotopic to $u_0$ and $\varphi_1$ is homotopic to $u_1$. This yields a contradiction to the $\varepsilon$-indecomposability of $\varphi$.
\end{proof}

\br \label{rmk:unif-distributed-energy-constants}
The proofs of \Cref{prop:unif-distributed-energy} and \Cref{prop:unif-distributed-energy-sphere} show that their corresponding constants $r_0$ depend only on $X,M,\varepsilon$ and $e_\mathbf{I}(\varphi)$, and that $r_0$ is decreasing with regard to the value of $e_\mathbf{I}(\varphi)$.
\er

\section{Proofs of main theorem and its consequences}\label{sec:proof-main-thm}
The goal of this section is to provide the proof of the main theorem and deduce its corollaries stated in \Cref{sec:Introduction}. Let $\mathbf{I}$ be a definition of energy as in \Cref{sec:energy-new}.
The following is a reformulation of our main result, \Cref{thm:main-introduction}, using the terminology established in the previous sections.

\bt\label{thm:decomposition-minimizers}
Let $M$ be a closed surface equipped with a Riemannian metric. Let $X$ be a compact quasiconvex metric space admitting a local quadratic isoperimetric inequality and such that every continuous map from $S^2$ to $X$ of sufficiently small diameter is null-homotopic. Then every continuous map $\varphi\colon M\to X$ has an iterated decomposition  such that $$e_\mathbf{I}(\varphi_0) + e_\mathbf{I}(\varphi_1) + \dots + e_\mathbf{I}(\varphi_k) = e_\mathbf{I}(\varphi)$$ and such that every $\varphi_i$ contains an $\mathbf{I}$-energy minimizer in its homotopy class.
\et

We need the following proposition in the proof of the theorem.

\bp\label{lem:energy-of-decomposition}
Let $M$ and $X$ be as in \Cref{thm:decomposition-minimizers}. Every iterated decomposition of a continuous map $\varphi\colon M\to X$ satisfies $$e_\mathbf{I}(\varphi)\leq e_\mathbf{I}(\varphi_0)+\dots+e_\mathbf{I}(\varphi_k).$$
\ep

\begin{proof} 
Let $\varphi_0\colon M\to X$ and $\varphi_1\colon S^2\to X$ be a decomposition of $\varphi$. Thus, there exists a disc $B\subset M$ such that $\varphi_0$ agrees with $\varphi$ on $M\setminus B$ and we can obtain $\varphi_1$ by gluing $\varphi_0|_B$ with $\varphi|_B$ along their common boundary.

Fix $\varepsilon>0$ and let $u_0\in\Lambda(\varphi_0)$ and $u_1\in\Lambda(\varphi_1)$ be such that $$E_\mathbf{I}(u_0)< e_\mathbf{I}(\varphi_0)+\varepsilon\quad\text{and}\quad E_\mathbf{I}(u_1)< e_\mathbf{I}(\varphi_1)+\varepsilon.$$ By modifying $u_0$ and $u_1$ in a way that adds an arbitrarily small amount of energy, we may assume that $u_0$ is constant in an open set $V\subset B$ and $u_1$ is constant on the homeomorphic image of $V$ in $S^2$.

Fix a point $p_0\in V$ and let $p_1$ be the corresponding point in $S^2$. Let $H_i$ be a free homotopy from $\varphi_i$ to $u_i$ for $i=0,1$. If we define $\gamma_0(t)\coloneqq H_0(p_0,1-t)$ and $\gamma_1(t)\coloneqq H_1(p_1,t)$ we can obtain a curve $\gamma^*$ that first goes along $\gamma_0$ and then $\gamma_1$, from $u_0(p_0)$ to $u_1(p_1)$. By quasiconvexity of $X$, we may find a Lipschitz curve $\gamma$ that is arbitrarily close to $\gamma^*$ and has the same endpoints. Then \Cref{lem:curves-small-diam-null-homotopic} implies that $\gamma$ is homotopic to $\gamma^*$ relative to its endpoints. 

Now let $V_0\subset V$ be a disc containing $p_0$, and choose $V_1\subset V_0$ to be a small disc such that $V_0\setminus V_1$ is conformally equivalent to an annulus whose two radii have a large ratio. We construct a continuous Sobolev map $u:M\to X$ as follows. Firstly, we define $u$ to agree with $u_0$ on $M\setminus V_0$. On the annulus $V_0\setminus V_1$, let $u$ coincide with a Lipschitz parametrization of $\gamma$. The energy of $u|_{V_0\setminus V_1}$ can be made to be smaller than $\varepsilon$ since we may choose the radius of $V_1$ to be arbitrarily small relative to the radius of $V_0$. On $V_1$, we define $u$ to agree with the map $u_1\circ \rho$, where $\rho:V_1\to S^2$ is a conformal map onto the complement of a small ball around $p_1$ in $S^2$. With this construction $u$ is a continuous Sobolev map freely homotopic to $\varphi$ and we may estimate its energy by
\begin{align*}
    E_\mathbf{I}(u) < E_\mathbf{I}(u_0)+E_\mathbf{I}(u_1)+\varepsilon <  e_\mathbf{I}(\varphi_0)+e_\mathbf{I}(\varphi_1)+3\varepsilon.
\end{align*}
Since $\varepsilon$ can be chosen arbitrarily small it follows that $e_\mathbf{I}(\varphi)\leq e_\mathbf{I}(\varphi_0)+e_\mathbf{I}(\varphi_1)$. From this, the statement in the proposition also holds for iterated decompositions.
\end{proof}

We use the proposition above together with the results established in the previous sections to prove our main result.

\begin{proof}[Proof of \Cref{thm:decomposition-minimizers}] 
Let $\varphi\colon M\to X$ be continuous and notice that $e_\mathbf{I}(\varphi)$ is finite due to \Cref{prop:Sobolev-homotopic-to-cont}. For $k\geq 0$ let $e_k$ be the infimal sum $e_\mathbf{I}(\varphi_0)+\dots+e_\mathbf{I}(\varphi_k)$ over $k$-step iterated decompositions of $\varphi$. If $\varphi$ does not have a $k$-step iterated decomposition then we set $e_k=\infty$. Notice that $e_0=e_\mathbf{I}(\varphi)$ and $e_{k+1}\geq e_k$ for all $k$ by \Cref{lem:energy-of-decomposition}. Furthermore, by \Cref{thm:small-area-spheres-contractible}  it holds that $e_k>e_\mathbf{I}(\varphi)$ for all $k$ large enough. Let $m$ be the biggest integer such that $e_m=e_\mathbf{I}(\varphi)$ and take a sequence of $m$-step iterated decompositions of $\varphi$ satisfying $$e_\mathbf{I}(\varphi_0^n)+\dots+e_\mathbf{I}(\varphi_m^n)\to e_\mathbf{I}(\varphi)$$ as $n$ tends to infinity. Set $\varepsilon\coloneqq2^{-1}\min\{\varepsilon_{\mathbf{I}},e_{m+1}-e_m\}>0$, where $\varepsilon_{\mathbf{I}}$ is defined as at the beginning of \Cref{sec:min-seq}. Clearly, $\varphi_i^n$ is $\varepsilon$-indecomposable for each $i$ and every sufficiently large $n$.

For all $i$ and $n$ choose $u_i^n\in \Lambda(\varphi_i^n)$ with $E_\mathbf{I}(u_i^n) \leq e_\mathbf{I}(\varphi_i^n) + {1}/{n}$. Then \Cref{prop:unif-distributed-energy}, \Cref{prop:unif-distributed-energy-sphere} and \Cref{rmk:unif-distributed-energy-constants} imply that by possibly precomposing the maps with a conformal diffeomorphism we may assume that there exists an $r_0>0$ such that, for sufficiently large $n$, the restriction of each $u_i^n$ to any ball of radius $r_0$ has $\mathbf{I}$-energy less than $\varepsilon_\mathbf{I}/5$.

After possibly passing to a subsequence we may assume, by the Rellich-Kondra\-chov compactness theorem \cite[Theorem 1.13]{KS93}, that each sequence $(u_i^{n})$ converges in $L^2$ to a Sobolev map $u_i$. By \Cref{thm:energy-minimizing-sequence}, each $u_i$ is continuous, satisfies $E_\mathbf{I}(u_i)=e_\mathbf{I}(u_i)$, and $u_i$ is homotopic to $u_i^{n}$  and thus to $\varphi_i^n$ for all sufficiently large $n$. Thus, for sufficiently large $n$, we have
\begin{align*}
    E_\mathbf{I}(u_0)+\dots +E_\mathbf{I}(u_m)&=e_\mathbf{I}(u_0)+\dots+e_\mathbf{I}(u_m)=e_\mathbf{I}(\varphi_0^{n})+\dots+e_\mathbf{I}(\varphi_m^{n}),
\end{align*}
which establishes the theorem.
\end{proof}

We end this section with the proofs of the corollaries to the main theorem stated in the introduction.

\begin{proof}[Proof of \Cref{cor:trivial-second-homotopy-group-intro} using \Cref{thm:decomposition-minimizers}]
Let $\varphi\colon M\to X$ be a continuous map. Since $\pi_2(X)$ is trivial, $\varphi$ does not have a decomposition, because this in particular would require the existence of an essential map from $S^2$. Then the only iterated decomposition of $\varphi$ is the zero-step iterated decomposition $\varphi_0=\varphi$ and thus \Cref{thm:decomposition-minimizers} guarantees the existence of an $\mathbf{I}$-energy minimizer in $\Lambda(\varphi)$.
\end{proof}

Before proving \Cref{thm:full-Sacks-Uhlenbeck-metric-space}, let us recall that there is a one-to-one correspondence between free homotopy classes of maps from $S^2$ to $X$ and orbits of the form $\pi_1(X)\cdot \alpha$ for $\alpha\in \pi_2(X)$, see \cite[Section 5]{SU-1981}.

\begin{proof}[Proof of \Cref{thm:full-Sacks-Uhlenbeck-metric-space}] 
Let $\Lambda$ be the family of free homotopy classes of continuous maps from $S^2$ to $X$ that are $\varepsilon$-indecomposable for some $\varepsilon>0$. Let $P$ be the subgroup of $\pi_2(X)$ generated by based representatives of  elements of $\Lambda$. 
By \Cref{thm:decomposition-minimizers}, every free homotopy class in $\Lambda$ contains an $\mathbf{I}$-energy minimizer, and hence it suffices to prove that $P=\pi_2(X)$. 

In order to prove this we argue by contradiction and assume that $\pi_2(X)\setminus P$ is non-empty. 
Let $e$ be the infimum of $e_{\mathbf{I}}(\varphi)$ over all $[\varphi]\in \pi_2(X)\setminus P$ and notice that $e\geq \varepsilon_{\mathbf{I}}>0$, where $\varepsilon_{\mathbf{I}}$ is as at the beginning of \Cref{sec:min-seq}. Fix $[\varphi]\in \pi_2(X)\setminus P$ such that $e_\mathbf{I}(\varphi) < e+\varepsilon_\mathbf{I}/2$. Since $\varphi$ is not $\varepsilon$-indecomposable for any $\varepsilon>0$ we can find a decomposition $\varphi_0\colon S^2\to X$ and $\varphi_1\colon S^2\to X$ of $\varphi$ such that
$$
e_\mathbf{I}(\varphi_0) + e_\mathbf{I}(\varphi_1)\leq e_\mathbf{I}(\varphi) + \varepsilon_\mathbf{I}/2< e+ \varepsilon_{\mathbf{I}}. 
$$
As $\varphi_0$ and $\varphi_1$ are essential we must have $e_\mathbf{I}(\varphi_i)\geq\varepsilon_\mathbf{I}$ for $i=0, 1$ and consequently $e_\mathbf{I}(\varphi_i)<e$; in particular, we have $\varphi_i\in \Lambda$. Thus, if $\varphi^*_i$ is a based representative of the free homotopy class of $\varphi_i$ then we have $\pi_1(X)[\varphi_i^*]\subset P$. Since $$\pi_1(X)[\varphi]\subset \pi_1(X)[\varphi_0^*]+\pi_1(X)[\varphi_1^*]\subset P,$$ compare with \cite[Section 5]{SU-1981}, we conclude that $[\varphi]\in P$, which is a contradiction.
\end{proof}

\def\cprime{$'$} \def\cprime{$'$} \def\cprime{$'$}


\begin{thebibliography}{10}

\bibitem{Amb90}
Luigi Ambrosio.
\newblock Metric space valued functions of bounded variation.
\newblock {\em Ann. Scuola Norm. Sup. Pisa Cl. Sci. (4)}, 17(3):439--478, 1990.

\bibitem{Avvakumov-Nabutovsky}
Sergey Avvakumov and Alexander Nabutovsky.
\newblock Boxing inequalities in banach spaces.
\newblock {\em preprint ArXiv: 2304.02709}, 2023.

\bibitem{Breiner-et-al}
Christine Breiner, Ailana Fraser, Lan-Hsuan Huang, Chikako Mese, Pam Sargent,
  and Yingying Zhang.
\newblock Existence of harmonic maps into {${\rm CAT}(1)$} spaces.
\newblock {\em Comm. Anal. Geom.}, 28(4):781--835, 2020.

\bibitem{Curtis-Fort-1957}
M.~L. Curtis and M.~K. Fort, Jr.
\newblock Homotopy groups of one-dimensional spaces.
\newblock {\em Proc. Amer. Math. Soc.}, 8:577--579, 1957.

\bibitem{EL78}
J.~Eells and L.~Lemaire.
\newblock A report on harmonic maps.
\newblock {\em Bull. London Math. Soc.}, 10(1):1--68, 1978.

\bibitem{EL88}
J.~Eells and L.~Lemaire.
\newblock Another report on harmonic maps.
\newblock {\em Bull. London Math. Soc.}, 20(5):385--524, 1988.

\bibitem{Futaki-1980}
Akito Futaki.
\newblock Nonexistence of minimizing harmonic maps from {$2$}-spheres.
\newblock {\em Proc. Japan Acad. Ser. A Math. Sci.}, 56(6):291--293, 1980.

\bibitem{Haj96}
Piotr Haj{\l}asz.
\newblock Sobolev spaces on an arbitrary metric space.
\newblock {\em Potential Anal.}, 5(4):403--415, 1996.

\bibitem{HL03}
Fengbo Hang and Fanghua Lin.
\newblock Topology of {S}obolev mappings. {II}.
\newblock {\em Acta Math.}, 191(1):55--107, 2003.

\bibitem{HKST15}
Juha Heinonen, Pekka Koskela, Nageswari Shanmugalingam, and Jeremy Tyson.
\newblock {\em Sobolev spaces on metric measure spaces}, volume~27 of {\em New
  Mathematical Monographs}.
\newblock Cambridge University Press, Cambridge, 2015.

\bibitem{Jost-1991-two-dim-var-problems}
J\"urgen Jost.
\newblock {\em Two-dimensional geometric variational problems}.
\newblock Pure and Applied Mathematics (New York). John Wiley \& Sons, Ltd.,
  Chichester, 1991.
\newblock A Wiley-Interscience Publication.

\bibitem{Jost94}
J\"{u}rgen Jost.
\newblock Equilibrium maps between metric spaces.
\newblock {\em Calc. Var. Partial Differential Equations}, 2(2):173--204, 1994.

\bibitem{Kar07}
M.~B. Karmanova.
\newblock Area and co-area formulas for mappings of the {S}obolev classes with
  values in a metric space.
\newblock {\em Sibirsk. Mat. Zh.}, 48(4):778--788, 2007.

\bibitem{Kir94}
Bernd Kirchheim.
\newblock Rectifiable metric spaces: local structure and regularity of the
  {H}ausdorff measure.
\newblock {\em Proc. Amer. Math. Soc.}, 121(1):113--123, 1994.

\bibitem{KS93}
Nicholas~J. Korevaar and Richard~M. Schoen.
\newblock Sobolev spaces and harmonic maps for metric space targets.
\newblock {\em Comm. Anal. Geom.}, 1(3-4):561--659, 1993.

\bibitem{Kuratowski1935}
Casimir Kuratowski.
\newblock Quelques problèmes concernant les espaces métriques
  non-séparables.
\newblock {\em Fundamenta Mathematicae}, 25(1):534--545, 1935.

\bibitem{Lem78}
Luc Lemaire.
\newblock Applications harmoniques de surfaces riemanniennes.
\newblock {\em J. Differential Geom.}, 13(1):51--78, 1978.

\bibitem{LS-2019}
Alexander Lytchak and Stephan Stadler.
\newblock Conformal deformations of {$\rm CAT(0)$} spaces.
\newblock {\em Math. Ann.}, 373(1-2):155--163, 2019.

\bibitem{LS-2020}
Alexander Lytchak and Stephan Stadler.
\newblock Improvements of upper curvature bounds.
\newblock {\em Trans. Amer. Math. Soc.}, 373(10):7153--7166, 2020.

\bibitem{LW16-harmonic}
Alexander Lytchak and Stefan Wenger.
\newblock Regularity of harmonic discs in spaces with quadratic isoperimetric
  inequality.
\newblock {\em Calc. Var. Partial Differential Equations}, 55(4):55:98, 2016.

\bibitem{LW15-Plateau}
Alexander Lytchak and Stefan Wenger.
\newblock Area minimizing discs in metric spaces.
\newblock {\em Arch. Ration. Mech. Anal.}, 223(3):1123--1182, 2017.

\bibitem{LW17-en-area}
Alexander Lytchak and Stefan Wenger.
\newblock Energy and area minimizers in metric spaces.
\newblock {\em Adv. Calc. Var.}, 10(4):407--421, 2017.

\bibitem{LW-intrinsic}
Alexander Lytchak and Stefan Wenger.
\newblock Intrinsic structure of minimal discs in metric spaces.
\newblock {\em Geom. Topol.}, 22(1):591--644, 2018.

\bibitem{LW-isoperimetric}
Alexander Lytchak and Stefan Wenger.
\newblock Isoperimetric characterization of upper curvature bounds.
\newblock {\em Acta Math.}, 221(1):159--202, 2018.

\bibitem{LW-param}
Alexander Lytchak and Stefan Wenger.
\newblock Canonical parameterizations of metric disks.
\newblock {\em Duke Math. J.}, 169(4):761--797, 2020.

\bibitem{LWY20}
Alexander Lytchak, Stefan Wenger, and Robert Young.
\newblock Dehn functions and {H}\"{o}lder extensions in asymptotic cones.
\newblock {\em J. Reine Angew. Math.}, 763:79--109, 2020.

\bibitem{Milnor65}
John~W. Milnor.
\newblock {\em Topology from the differentiable viewpoint}.
\newblock University Press of Virginia, Charlottesville, VA, 1965.
\newblock Based on notes by David W. Weaver.

\bibitem{Res97}
Yu.~G. Reshetnyak.
\newblock Sobolev classes of functions with values in a metric space.
\newblock {\em Sibirsk. Mat. Zh.}, 38(3):657--675, iii--iv, 1997.

\bibitem{SU-1981}
J.~Sacks and K.~Uhlenbeck.
\newblock The existence of minimal immersions of {$2$}-spheres.
\newblock {\em Ann. of Math. (2)}, 113(1):1--24, 1981.

\bibitem{SY79}
R.~Schoen and Shing~Tung Yau.
\newblock Existence of incompressible minimal surfaces and the topology of
  three-dimensional manifolds with nonnegative scalar curvature.
\newblock {\em Ann. of Math. (2)}, 110(1):127--142, 1979.

\bibitem{Soultanis-Wenger-2022}
Elefterios Soultanis and Stefan Wenger.
\newblock Area minimizing surfaces in homotopy classes in metric spaces.
\newblock {\em Trans. Amer. Math. Soc.}, 375(7):4711--4739, 2022.

\bibitem{White86}
Brian White.
\newblock Infima of energy functionals in homotopy classes of mappings.
\newblock {\em J. Differential Geom.}, 23(2):127--142, 1986.

\bibitem{White88}
Brian White.
\newblock Homotopy classes in {S}obolev spaces and the existence of energy
  minimizing maps.
\newblock {\em Acta Math.}, 160(1-2):1--17, 1988.

\end{thebibliography}
\end{document}